\documentclass[a4paper,11pt]{article}
\usepackage[english]{babel}
\usepackage[latin1]{inputenc}
\usepackage{graphicx}
\usepackage{amscd}
\usepackage{amsfonts}
\usepackage{latexsym}
\usepackage{amsthm}
\usepackage{amssymb,amsmath}
\usepackage{enumerate}
\usepackage{verbatim}

\usepackage{import}

\usepackage[usenames]{color}

\usepackage{subfigure}

\title{\bf On a Family of Rational Perturbations of the Doubling Map}

\author{\small Jordi Canela \thanks{ The three authors were supported by MTM2011-26995-C02-02, CIRIT 2009-SGR792 and the EU project MRTN-CT-2006-035651. The first author was also supported by the Spanish government grant FPU AP2009-4564} \\
{\small Dep.~de Matem\`atica Aplicada i An\`alisi}\\
{\small Universitat de Barcelona}\\
{\small Gran Via de les Corts Catalanes, 585}\\
{\small 08005 Barcelona, Spain} \and 
{\small N\'uria Fagella $^*$}\\
{\small Dep.~de Matem\`atica Aplicada i An\`alisi}\\
{\small Universitat de Barcelona}\\
{\small Gran Via de les Corts Catalanes, 585}\\
{\small 08005 Barcelona, Spain} \and 
{\small Antonio Garijo $^*$}\\
{\small Dept. d'Enginyeria Inform\`atica i Matem\`atiques}\\
{\small Universitat Rovira i Virgili}\\
{\small Av. Pa\"isos Catalans 26}\\
{\small Tarragona 43007, Spain}  }

\newtheorem{teor}{Theorem} [section]
\newtheorem{propo}[teor]{Proposition}

\newtheorem{lemma}[teor]{Lemma}
\newtheorem{co}[teor]{Corollary}
\newtheorem{conj}[teor]{Conjecture}
\newtheorem*{teoremB}{Theorem B}
\newtheorem*{teoremA}{Theorem A}
\newtheorem*{teoremC}{Theorem C}

\theoremstyle{definition}
\newtheorem{defin}[teor]{Definition}

\newsavebox{\savepar}

\newcommand{\com}{\mathbb{C}}
\newcommand{\real}{\mathbb{R}}

\newcommand{\nat}{\mathbb{N}}

\newcommand{\dis}{\mathbb{D}}
\newcommand{\cercle}{\mathbb{S}^1}
\newcommand{\re}{\rm{Re}}
\newcommand{\im}{\rm{Im}}

\newcommand{\modul}{\rm{mod}\;}

\date{\today}

    {\end{list}}

\begin{document}
\maketitle

\begin{abstract}
 {\noindent \small The goal of this paper is to investigate the parameter plane of a rational family of perturbations of the doubling map given by the Blaschke products  $B_a(z)=z^3\frac{z-a}{1-\bar{a}z}$. First we study the basic properties of these maps such as the connectivity of the Julia set as a function of the parameter $a$. We use techniques of quasiconformal surgery to explore the relation between certain members of the family and the degree 4 polynomials $\left(\overline{\overline{z}^2+c}\right)^2+c$. In parameter space, we classify the different hyperbolic components according to the critical orbits and we show how to parametrize those of disjoint type.}

\vspace{0.5cm}
\textit{Keywords: holomorphic dynamics, Blaschke products, circle maps, polynomial-like mappings.}
\end{abstract}


\section{Introduction}

 Given a rational map $f:\widehat{\com}\rightarrow \widehat{\com}$, where $\widehat{\com}=\com\cup \{\infty\}$ denotes the Riemann sphere, we consider the dynamical system given by the iterates of $f$. The Riemann sphere splits into two totally $f$-invariant subsets: the  \textit{Fatou set} $\mathcal{F}(f)$, which  is defined to be the set of points $z\in\widehat{\com}$ where the family $\{f^n, n\in\mathbb{N}\}$ is normal in some neighbourhood of $z$, and  its complement,  the \textit{Julia set} $\mathcal{J}(f)$. The dynamics of the points in $\mathcal{F}(f)$ are stable whereas the dynamics in $\mathcal{J}(f)$ present chaotic behaviour. The Fatou set $\mathcal{F}(f)$ is  open and therefore $\mathcal{J}(f)$ is  closed. Moreover, if the degree of the rational map $f$ is greater or equal than $2$, then the Julia set $\mathcal{J}(f)$ is not empty and is the closure of the set of  repelling fixed points of $f$.
 
  The connected components of $\mathcal{F}(f)$, called \textit{Fatou components}, are mapped under $f$ among themselves. D.~Sullivan \cite{Su} proved that any Fatou component of a rational map is either periodic or preperiodic. By means of the Classification Theorem (see e.g.\ \cite{Mi1}), any periodic Fatou component of a rational map  is either the basin of attraction of an attracting or parabolic cycle or is a simply connected rotation domain (a Siegel disk) or is a doubly connected rotation domain (a Herman Ring). Moreover, any such component is somehow related to a \textit{critical point}, i.e.\ a point $z\in\widehat{\com}$ such that $f'(z)=0$. Indeed, the basin of attraction of an attracting or parabolic cycle contains, at least, a critical point whereas Siegel disks and Herman rings have critical orbits accumulating on their boundaries. For a background on the dynamics of rational maps we refer to \cite{Mi1} and \cite{Bear}.

The aim of this paper is to study the dynamics of the degree 4 Blaschke products given by 

\begin{equation}\label{blasformula}
B_{a}(z)=z^3\frac{z-a}{1-\bar{a}z},
\end{equation}

\noindent where $a,z\in\com$.  This Blaschke family restricted to $\cercle$ is the rational analogue of the double standard family $A_{\alpha,\beta}(\theta) =2\theta+\alpha+(\beta/\pi)\sin(2\pi\theta) \hspace{0.2cm} (\mbox{mod} \hspace{0.2cm} 1)$ of periodic perturbations of the doubling map on $\cercle$. Indeed, when $|a|$ tends to infinity, the products $B_a$ tend to $e^{4\pi i t}z^2$ uniformly on compact sets of  the punctured plane $\com^*=\com\setminus\{0\}$, where $t\in \real/\mathbb{Z}$ denotes the argument of $a$. Double standard maps extend to entire transcendental self-maps of $\com^*$. Although there is no explicit simple expression for the restriction of $B_a$ to $\cercle$, the global dynamics are simpler than in the transcendental case. The double standard family has been studied in several papers such as  \cite{MiRo1,MiRo2}, \cite{De} and \cite{dLSS}.

For all values of $a\in\com$, the points $z=0$ and $z=\infty$ are superattracting fixed points of local degree 3. We denote by $A(0)$ and $A(\infty)$ their basins of attraction and by $A^*(0)$ and $A^*(\infty)$ their immediate basins of attraction, i.e.\ the connected components of the basins containing the superattracting fixed points. If $|a|\leq 1$, $A(0)=A^*(0)=\dis$ and $A(\infty)=A^*(\infty)=\widehat{\com}\setminus\overline{\dis}$ and hence $\mathcal{J}(B_a)=\cercle$ (see Lemma \ref{juliacercle}). If $|a|>1$, there are two critical points $c_{\pm}$. If $|a|=2$, the two critical points collide in a single one. If $|a|>2$ the two critical points are symmetric with respect to $\cercle$ (see Section \ref{introblas}) and so are their orbits.  Consequently,  if $|a|>2$ it is enough to control one of the critical orbits in order to understand the possible dynamics of $B_a$. On the other hand, if $1<|a|<2$ the Blaschke product $B_a$ has two different free critical points contained in $\cercle$ which may lead to different dynamics.

The connectivity of the Julia set is a focus of attention when studying a family of rational functions (see e.g.\ \cite{Shi}, \cite{Pr}, \cite{Pi} and \cite{DR}). It is known that, given a polynomial $P$, its Julia set $\mathcal{J}(P)$ is connected if and only if it has no free critical point captured by the basin of attraction of infinity (see \cite{Mi1}). However, such a  classification does not exist for general rational maps which, unlike polynomials, may have Herman Rings  and even cantor sets of Jordan curves \cite{McM}. It turns out that the family $B_a$ shares some of the features of polynomials in this respect such as the non existence of Herman rings (see Proposition $\ref{noHR}$). We also prove the following criterium.

\begin{teoremA}

Given a Blaschke product $B_a$ as in $(\ref{blasformula})$, the following statements hold:
\begin{enumerate}[(a)]
\item If $|a|\leq 1$, then $\mathcal{J}(B_a)=\cercle$.
\item If $|a|>1$, then $A(\infty)$ and $A(0)$ are simply connected if and only if $c_+\notin A^*(\infty)$.
\item If $|a|\geq 2$, then every Fatou component $U$ such that $U\cap A(\infty)=\emptyset$ and $U\cap A(0)=\emptyset$ is simply connected.
\end{enumerate}
 Consequently, if $|a|\geq 2$, then $\mathcal{J}(B_a)$ is connected if and only if $c_+\notin A^*(\infty)$.

\end{teoremA}

 Next, we focus on the basins of attraction of attracting or parabolic cycles not contained in $\cercle$, other than $0$ or $\infty$. These may only exist when $|a|>2$. We distinguish two cases depending on the location of the cycles with respect to $\dis$.
 
 On the one hand, if $B_a$, $|a|>2$,  is such that  has no attracting or parabolic cycle in $\cercle$, we can relate the Blaschke product with a cubic polynomial preserving the dynamics of all orbits contained in $\com\setminus \overline{\dis}$. Indeed, for such parameters,  $B_{a}|_{\cercle}$ is quasisymmetrycally conjugate to the doubling map $\theta\rightarrow2\theta \;(\modul 1)$ and  a quasiconformal surgery, consisting in gluing a superattracting cycle in $\dis$,  can be performed  obtaining cubic polynomials of the form $M_b(z)=bz^2(z-1)$ with $b\in\com$ (see \cite{Pe} and Section \ref{surgery}). These polynomials have been the object of study of several papers (see e.g.\ \cite{Mi2} and \cite{Ro}).
   This surgery establishes a conformal conjugacy between $M_b$ and $B_a$ on the set of points which never enter  $\dis$ under iteration of $B_a$ and the points which are not attracted to $z=0$ under iteration of $M_b$. In particular, if $B_a$ has an attracting or parabolic cycle contained in $\com\setminus\overline{\dis}$, this surgery conjugates $B_a$ with $M_b$ conformally in its basin of attraction.

On the other hand, if $B_a$ has a periodic cycle with points both inside and outside $\dis$ the situation is different. Although the previous surgery construction is still possible, a lot of information is lost since, under the new map, the critical point is always captured by the basin of $z=0$. Parameters for which the orbit of  $c_+\in\com\setminus\overline{\dis}$ enters the unit disk at least once are called \textit{swapping parameters} and connected components of the set of swapping parameters are called \textit{swapping regions}. Inside these regions, the non holomorphic dependence of $B_a$ on  the parameter $a$  gives rise to what appear to be small copies of the Tricorn, the bifurcation locus of the antiholomorphic family (i.e.\ holomorphic on the variable $\overline{z}$) of polynomials $p_c(z)=\overline{z}^2+c$ (see \cite{Cr} and Figure \ref{swapregion} (a)).
J.~Milnor \cite{Mi3} showed that a similar situation takes place for real cubic  polynomials introducing the concept of antipolynomial-like mapping. We distinguish between two types of attracting cycles for swapping parameters. We say that a parameter is \textit{bitransitive} if it has a cycle whose basin of attraction contains the two free critical points.  We say that a parameter is \textit{disjoint} if there are two different cycles other than zero or infinity.
The very special dynamics taking place for swapping parameters allow us to build a polynomial-like mapping in a neighbourhood of every bitransitive or disjoint swapping parameter. A polynomial-like map is a locally defined map which is conjugate to a polynomial (see Definition \ref{defpollike}).

\begin{teoremB}
Let $a_0$ be a swapping parameter with an attracting or parabolic cycle of period $p>1$. Then, there is an open set $W$  containing $a_0$ and $p_0>1$ dividing $p$  such that, for every $a\in W$, there exist two open sets $U$ and $V$ with $c_+\in U$  such that $\left(B_a^{p_0}; U, V\right)$ is a polynomial-like map. Moreover,
\begin{enumerate}[(a)]
\item If $a_0$ is bitransitive,  $\left(B_a^{p_0}; U, V\right)$ is hybrid equivalent to a polynomial of the form $p^2_c(z)=\left(z^2+\overline{c}\right)^2+c$.
\item If $a_0$ is disjoint, $\left(B_a^{p_0}; U, V\right)$ is hybrid equivalent to a polynomial of the form $p^2_c(z)=\left(z^2+\overline{c}\right)^2+c$ or of the form $z^2+c$.
\end{enumerate}

\end{teoremB}

It is known that the boundary of every bounded Fatou component of a polynomial, with the exception of Siegel disks, is a Jordan curve \cite{RoYi}. This is not true however for rational functions although it can be stablished under some conditions (e.g.\ postcritically finite among others \cite{Pi}). In our case, as a consequence of the two previous constructions, we know that the boundary of every connected component of the basin of attraction of an attracting or parabolic cycle of $B_a$ not contained in $\cercle$ and other than $z=0$ and $z=\infty$ is a Jordan curve (see Proposition \ref{jordancurveattr}). Indeed, if $B_a$ has such a cycle, the previous constructions provide a conjugation between $B_a$ and a polynomial which sends the immediate basin of attraction of the cycle of $B_a$ to the immediate basin of attraction of a bounded cycle of the polynomial.  

 A rational map is \textit{hyperbolic} if all its critical points are attracted to attracting cycles. A hyperbolic component is a connected component of the open set $\mathcal{H}=\{a| B_a \mbox{ is hyperbolic} \}$. The parametrization of hyperbolic components of rational functions which depend holomorphically on their parameters is well known (see \cite{DH2} and \cite{BF}). If the family of functions does not depend holomorphically on parameters, some extra difficulties appear. S.~Nakane and D.~Schleicher \cite{NaSh} study the parametrization of hyperbolic components with cycles of even period for the family of antipolynomials  $p_{c,d}(z)=\overline{z}^d+c$. We focus on  the parametrization of hyperbolic components with disjoint parameters using different methods than the ones of \cite{NaSh}. Notice that, due to the symmetry of $B_a$, disjoint cycles are symmetric with respect to $\cercle$ and therefore have the same period and conjugate multiplier (see Theorem \ref{clashyp}). Hence, given a hyperbolic component $U$ with disjoint parameters,  it makes sense to define the multiplier map $\Lambda:U\rightarrow \dis$ as the map which sends every $a\in U$ to the multiplier of the attracting cycle whose basin captures the critical orbit of $c_+$.

\begin{teoremC}
Let $U$ be a disjoint hyperbolic component such that $U\subset \{a\in\com; |a|>2\}$. Then, the multiplier map is a homeomorphism between $U$ and the unit disk.

\end{teoremC}

Since the multiplier of any bitransitive cycle is a non-negative real number (see Proposition \ref{bitransmult}), the previous result does not hold for bitransitive components. This phenomena had already been noticed in $\cite{NaSh}$ for the polynomials $p_{c,d}^2$.

In Section \ref{stools} we introduce some notation and useful results to prove Theorems A, B and C. In Section \ref{introblas} we describe the basic properties of the Blaschke family $B_a$ and prove Theorem~A. In Section~\ref{parameterplane} we study the parameter plane of the family: we describe the different kinds of hyperbolic dynamics that may occur depending on the behaviour of the free critical points, we introduce the relation with the family of cubic polynomials $M_b$, describe the dynamics that can take place along the swapping regions proving Theorem~B and finally we prove Theorem~C.

\section{Preliminaries and tools}\label{stools}

Given a rational map $f$, we denote by $<z_0>:=\{z_0,\cdots,f^{p-1}(z_0)=z_{p-1}\}$  an attracting or parabolic cycle of period $p\geq1$, where $f(z_i)=z_{i+1}$ with subindeces taken modulus $p$. We denote with $A(<z_0>)$ the basin of attraction of the cycle whereas  $A^*(<z_0>)$ denotes  its immediate basin of attraction, i.e.\ the connected components of $A(<z_0>)$ which contain a point of $<z_0>$. With $A^*(z_q)$ we denote the connected component of $A^*(<z_0>)$ containing $z_q$. The marked point $z_0$ of the cycle is usually taken so that $A^*(z_0)$ contains a critical point.

\subsubsection*{The Riemann-Hurwitz formula}
 When dealing with the simple connectivity of open sets, it is useful to consider the Riemann-Hurwitz formula (see \cite{Bear}). It can be stated as follows. 

\begin{teor}[Riemann-Hurwitz Formula]\label{riemannhurwitz}
Let $U$ and $V$ be two connected domains  of $\widehat{\com}$ of finite connectivity $m_{U}$ and $m_{V}$ and let $f:U\rightarrow V$ be a degree $k$ proper map branched over $r$ critical points counted with multiplicity. Then
$$m_{U}-2=k(m_{V}-2)+r.$$
\end{teor}

The following corollary is used several times along the paper.

\begin{co}\label{connpreim}
Let $f$ be a rational map and let $V$ be a simply connected domain. Let $U$ be a connected component of $f^{-1}(V)$. If $U$ contains at most one critical point (of arbitrary multiplicity), then $U$ is simply connected.
\end{co}
\proof
By construction, $f|_{U}:U\rightarrow V$ is proper. Let $r$ be the multiplicity of the critical point. Then, $f|_{U}$ has at least degree $r+1$. By The Riemann-Hurwitz formula, since $m_V=1$,  we have $m_U-2\leq -(r+1)+r=-1$. Since $m_U$ is at least 1, we conclude that it is indeed 1 and $U$ is simply connected.
\endproof

\subsubsection*{Polynomial and antipolynomial-like mappings}\label{santipol}
 
 The theory of polynomial like mappings \cite{DH1}, introduced by A.~Douady and J.~Hubbard,  allows us to understand why copies of the filled Julia set of polynomials appear in the dynamical planes of rational maps or even entire transcendental or meromorphic maps.  

\begin{defin}\label{defpollike} A triple $(f; U, V)$ is called a \textit{polynomial-like} (resp.\ \textit{antipolynomial-like}) mapping of degree $d$ if $U$ and $V$ are bounded simply connected subsets of the plane, $\overline{U}\subset V$ and $f:U\rightarrow V$ is holomorphic (resp.\ antiholomorphic) and proper of degree $d$. Its  \textit{filled Julia set} is defined as
 $$\mathcal{K}_f=\bigcap_{n>0}f^{-n}(V)=\{z\in U |\,f^n(z)\in U \hspace{0.2cm} \forall n\geq 0\}.$$
\end{defin} 

Observe that given any polynomial (resp.\ antipolynomial) $P$ of degree $d$, there exists a disk $\dis_R$ of radius $R>0$ so that $(P; \dis_R, P(\dis_R))$ is a polynomial-like mapping (resp.\ antipolynomial-like mapping).

\begin{defin}
We say that two (anti)polynomial-like maps $(f; U, V)$  and $(f'; U', V')$ are hybrid equivalent if  there exist neighbourhoods $U_f$ and $U_{f'}$ of $\mathcal{K}_f$ and $\mathcal{K}_{f'}$, respectively, and a quasiconformal conjugation $\phi: U_f\rightarrow U_{f'}$ between $f$ and $f'$ such that $\overline{\partial} \phi=0$ almost everywhere in $\mathcal{K}_f$. 
\end{defin}

Polynomial-like mappings, as their name indicates, behave locally as polynomials do. This is the content of the  Straightening Theorem \cite{DH1}.

\begin{teor}[The Straightening Theorem]\label{straightening}
Every polynomial-like mapping $(f; U, V)$ of degree $d$ is hybrid equivalent to a polynomial $P$ of degree $d$. If $\mathcal{K}_f$ is connected, $P$ is unique up to affine conjugation.

\end{teor}
 
 The antipolynomial-like theory was first introduced by Milnor (see \cite{Mi2}) in order to study why small copies of the Tricorn appear in the parameter plane of real cubic polynomials. Hubbard and Schleicher (see \cite{HS}) used this theory afterwards in the study of the Multicorns, the parameter plane of the antipolynomial maps $p_{c,d}(z)=\overline{z}^d+c$. They stated the Antiholomorphic Straightening Theorem.

\begin{teor}[The Antiholomorphic Straightening Theorem]\label{antistraightening}
Every antipol\-ynomial-like mapping $(f; U, V)$ of degree $d$ is hybrid equivalent to an antipolynomial $P$ of degree $d$. If $\mathcal{K}_f$ is connected, then $P$ is unique up to affine conjugation.

\end{teor}

 \subsubsection*{Conjugation with the doubling map}

Along the paper it is important to know whether the Blaschke products $B_a$ are conjugate to the doubling map of the circle $\theta\rightarrow 2\theta\; (\mbox{mod} \hspace{0.2cm} 1)$, where $\theta \in \real/\mathbb{Z}$ (equivalently given by $z\rightarrow z^2$, where $z\in\cercle$). The following result tells us that the lift $H$ of an increasing covering map $f$ of the circle of degree $2$ is semiconjugate  to the doubling map $\theta\rightarrow 2\theta$. See \cite{deMVa, MiRo1, De} for further details.

\begin{lemma}\label{semiconj}
Let $F:\real\rightarrow\real$ be a continuous  increasing map. Suppose that $F(x+k)=F(x)+2k$ for any integer $k$  and for any real $x$. Then, the limit
$$H(x)=\lim_{n\rightarrow\infty}\frac{F^n(x)}{2^n}$$
exists uniformly on $x$.  This map $H$ is  increasing, continuous and satisfies $H(x+k)=H(x)+k$ for any integer $k$ and for any real $x$.  $H$ semiconjugates $F$ with the multiplication by $2$, i.e.\ $H(F(x))=2H(x)$ for any real $x$. Moreover, the map $H$ is unique up to constant $k\in\nat$ and sends points of period d to points of period d.
\end{lemma}

If $|a|\geq2$, the circle map $B_{a}|_{\cercle}$ is a degree 2 cover. We denote by $H_a$ the lift of the map $h_a:\cercle \rightarrow \cercle$ which semiconjugates $B_{a}|_{\cercle}$ with the doubling map in $\cercle$. Since $H_a$ is unique up to constant $k\in\nat$, the map $h_a$ is unique.

 \begin{defin}
 An orientation preserving map $h:\mathbb{S}^1\rightarrow \mathbb{S}^1$ is \textit{quasisymmetric} if $h$ is injective and such that, for $z_1,z_2,z_3\in \mathbb{S}^1$,
 $$\mbox{if} \hspace{0.2cm} |z_1-z_2|=|z_2-z_3|\Rightarrow\frac{1}{M}\leq\frac{|h(z_1)-h(z_2)|}{|h(z_2)-h(z_3|)}\leq M$$
 for some $M>0$.
 
\end{defin}

It is not difficult to prove that, if $|a|>3$, $B_{a}|_{\cercle}$ is and expanding map and, therefore, $h_a$ is a quasisymmetric conjugacy (cf.\ \cite{deMVa}). The next theorem,  due to C.~Petersen \cite{Pe}, gives some conditions weaker than expansivity which guarantee that $h_a$ is a quasisymmetric conjugacy. Recall that the $\omega$-limit set $\omega(z)$ of a point $z\in\com$ is defined to be the accumulation set of the orbit of $z$. A point $z$ is called recurrent if and only if $z\in\omega(z)$.

\begin{teor}\label{conjpe}
Let $B:\widehat{\com}\rightarrow\widehat{\com}$ be a Blaschke product with poles in $\dis$ such that the restriction $B:\cercle \rightarrow\cercle$ is a (positively oriented) degree $d\geq 2$ covering and such that $\cercle$ contains no non repelling periodic point and $\omega(c)\cap\cercle=\emptyset$ for every recurrent critical point c.
Then,  $B:\cercle \rightarrow\cercle$ is quasisymmetrically conjugate to $R_d(z)=z^d$.
\end{teor}

\section{Dynamical plane of the Blaschke family}\label{introblas}

We consider the degree 4 Blaschke products of the form 
\begin{equation}\label{blasformula2}
 B_{a,t}(z)=e^{2\pi it}z^3\frac{z-a}{1-\bar{a}z},
 \end{equation}

\noindent where $a \in \com$ and $t \in \mathbb{R}/\mathbb{Z}$. As all Blaschke products, the family $B_{a,t}$ leaves $\cercle$ invariant and therefore its members are symmetric with respect to the unit circle, i.e.\ $B_{a,t}(z)=\mathcal{I} \circ B_{a,t}\circ \mathcal{I}(z)$ where $\mathcal{I}(z)=1/\bar{z}$.

 The next lemma tells us that, for the purpose of classification,  we can get rid of the parameter $t$. The proof is straightforward.

\begin{lemma}\label{conjblas}
Let $\alpha\in \real$ and let $\eta(z)=e^{-2\pi i\alpha}z$. Then $\eta$ conjugates the maps $B_{a,t}$ and $B_{a e^{-2\pi i\alpha}, t+3\alpha}$. In particular we have that $B_{a,t}$ is conjugate to $B_{a e^{\frac{2\pi i t}{3}}, 0}$.

\end{lemma}

Hence, we focus on the study of the family $B_{a}(z)=z^3\frac{z-a}{1-\bar{a}z}$  (\ref{blasformula}) for values $a,z\in\com$.

  Given the fact that these rational maps have degree $4$, there are $6$ critical points counted with multiplicity. The fixed points $z=0$ and $z=\infty$ are critical points of multiplicity 2 and hence superattracting fixed points of local degree 3. The other two  critical points, denoted by $c_{\pm}$, are given by 

\begin{equation}\label{criticalpoints}
c_{\pm}:=c_{\pm}(a):=a \cdot \frac{1}{3|a|^2}\left(2+|a|^2\pm\sqrt{(|a|^2-4)(|a|^2-1)}\right).
\end{equation}

The critical points $c_+$ and $c_-$ are free and satisfy $|c_+|\geq 1$ and $|c_-|\leq 1$. If they are not in $\cercle$ (i.e.\ $|a|>2$ or $|a|<1$),  they and their orbits are symmetric with respect to $\cercle$.  The following result proves that the critical points determine the parameter if $|a|\geq 2$ or $|a|<1$.

\begin{lemma}\label{continuitya}
Given a Blaschke product $B_{a,t}$ as in (\ref{blasformula2}) with $|a|\geq2$ or $|a|< 1$, the parameter $a$ is continuously determined by the critical points $c_{\pm}$. Moreover, if the image $B_a(z_0)\neq\{0,\infty\}$ of a point $z_0$ is fixed, the parameter $t$ is continuously determined by $a$.

\end{lemma}
\proof
The continuous dependence of $t$ with respect to $a$ is clear. Let $a=r_a e^{2\pi i \alpha}$, where $\alpha\in \mathbb{R}/\mathbb{Z}$ and  $r_a\geq 2$ (resp.\ $r_a<1 $). It follows from (\ref{criticalpoints}) that the critical points  $c_+$ and $c_-$ have the same argument $\alpha$ as $a$. It is left to see that $r_a$ depends continuously on $|c_+|=r_c$. It follows from symmetry that $|c_-|=1/r_c$. Consider $R(r_c)=r_c+1/r_c$. For $r_c\geq 1$, $R$ is a strictly increasing function which satisfies $R(1)=2$. Using  (\ref{criticalpoints}) we have:
$$R(r_c)e^{2\pi i \alpha}=c_+ + c_-=\frac{2a}{3|a|^2}(2+|a|^2)=\frac{2}{3}\frac{r_a e^{2\pi i \alpha}}{r_a^2}(2+r_a^2),$$

\noindent and, therefore, $r_a\cdot R(r_c)=2(2+r_a^2)/3$. This quadratic equation yields two solutions $r_{a_{\pm}}=(3R\pm\sqrt{9R^2-32})/4$. The solution $r_{a_+}(R)$ takes the value $2$ for $R=2$ and is strictly increasing and tending to infinity when $R$ tends to infinity. The solution $r_{a_-}(R)$ takes the value $1$ for $R=2$ and is strictly decreasing and tending to zero when $R$ tends to infinity. Therefore, each critical point $c_+\in \com$, $|c_+|\geq1$ (resp.\ $|c_+|>1)$, determines continuously a unique parameter $a$ such that $|a|\geq 2$ (resp.\ $|a|<1$).

\endproof

Another relevant point to discuss is what type of dynamics may occur in $\cercle$. It follows directly from the invariance of $\cercle$ under $B_a$ and the fact that $B_a:\cercle\rightarrow\cercle$ cannot be a degree $1$ covering (and hence conjugate to an irrational rotation) that any point $z\in\cercle$ either belongs to the Julia set or is attracted to an attracting or a parabolic periodic orbit $\{z_0,...,z_{p-1}\}\in\cercle$.

\begin{figure}[hbt!]
    \centering
    \subfigure[\scriptsize{Case $|a|<1$} ]{
   			 \def\svgwidth{170pt}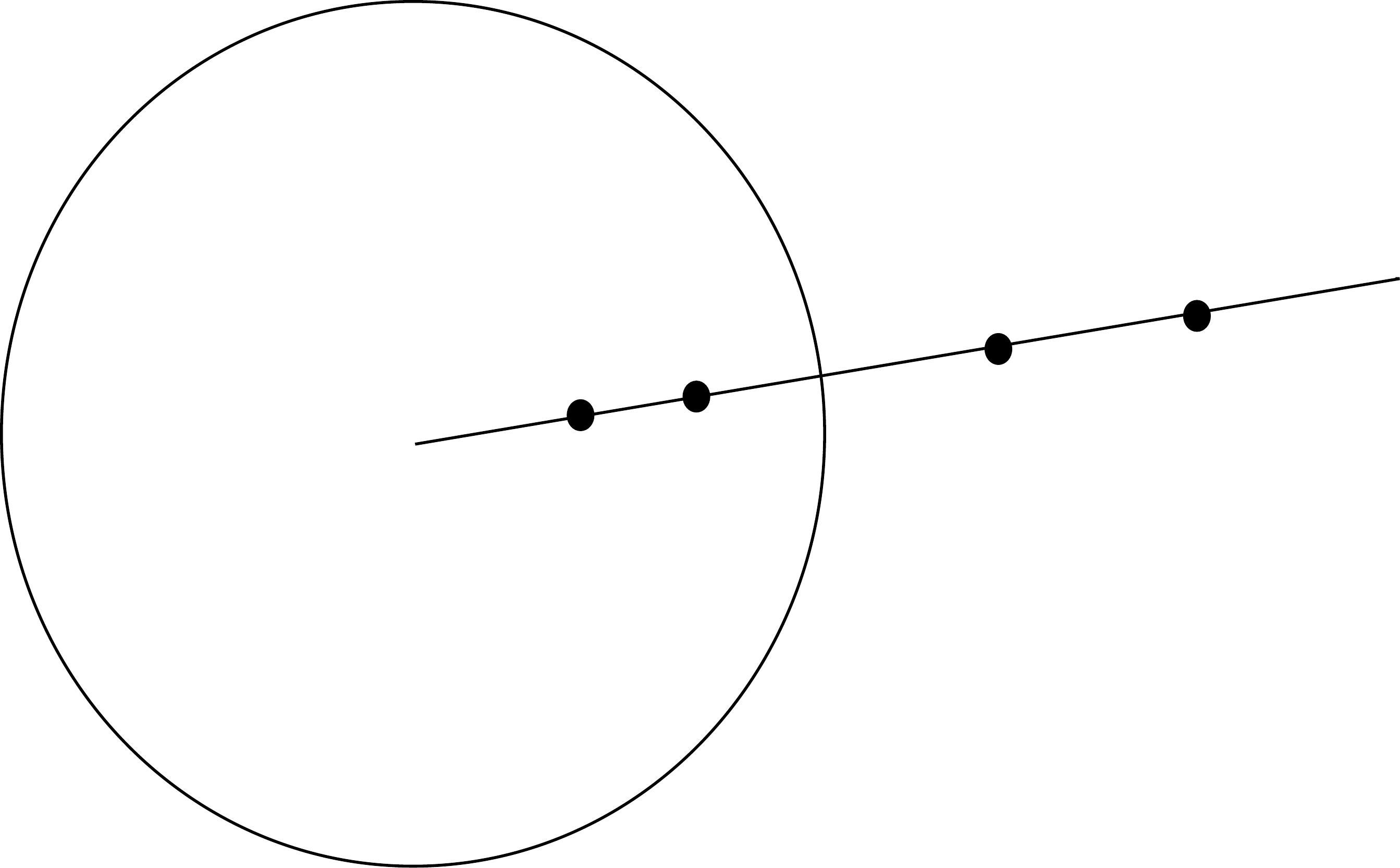
     }
    \subfigure[\scriptsize{Case $1<|a|<2$} ]{
    		 \def\svgwidth{170pt}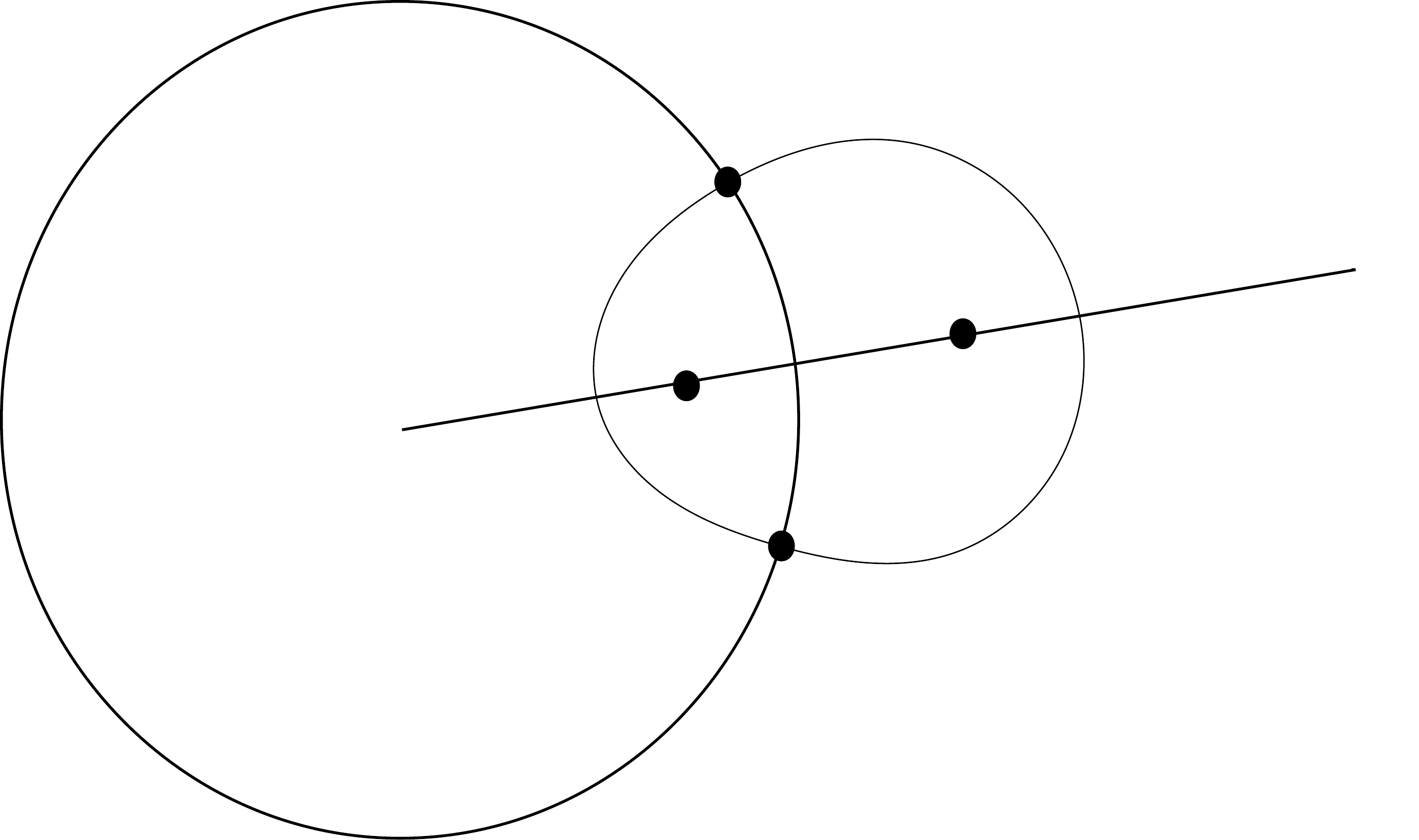}
    \hspace{0.1in}
    \subfigure[\scriptsize{Case $|a|=2$} ]{
 		  \def\svgwidth{170pt}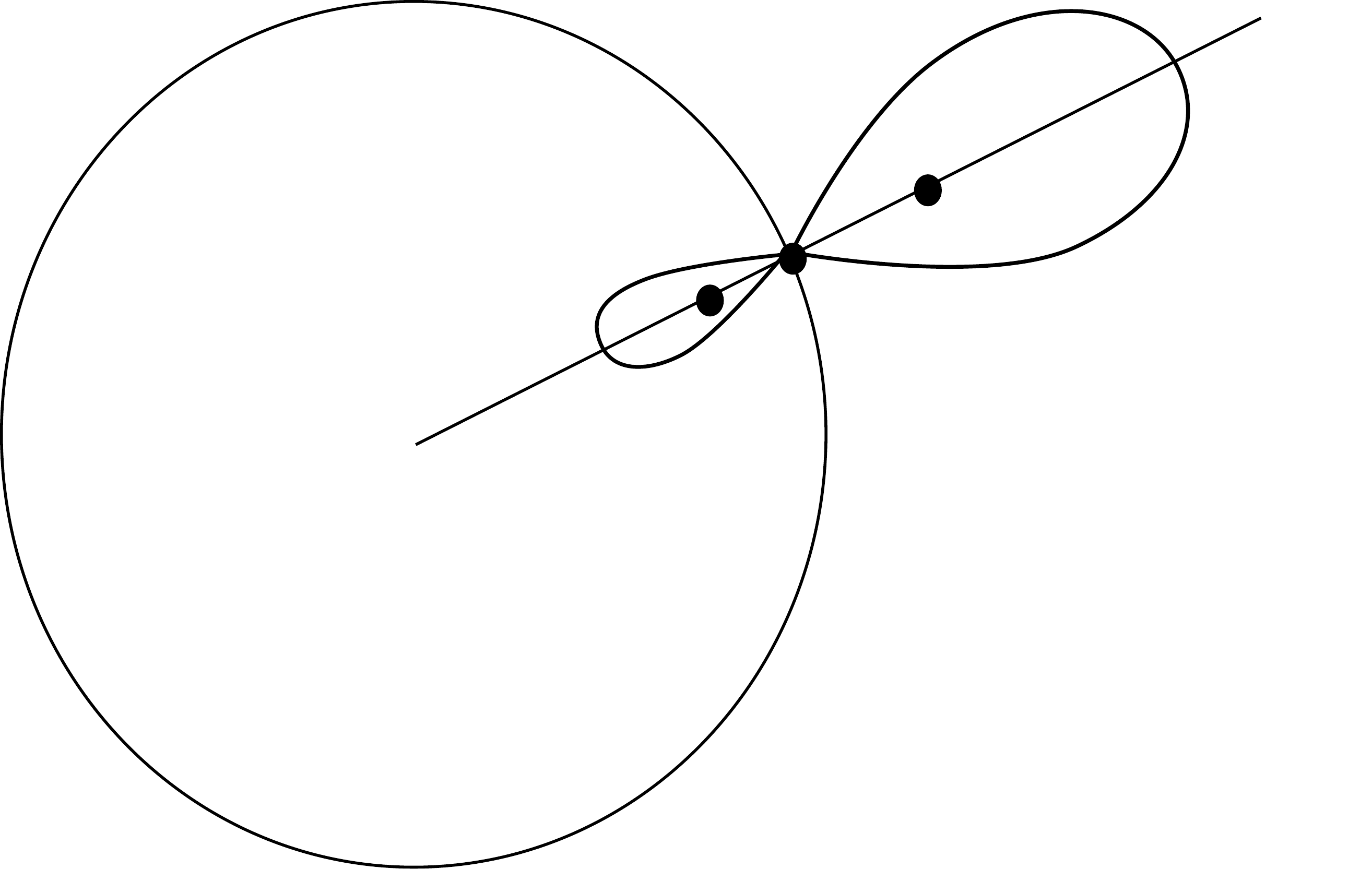}
    \subfigure[\scriptsize{Case $|a|>2$} ]{
   		   \def\svgwidth{170pt}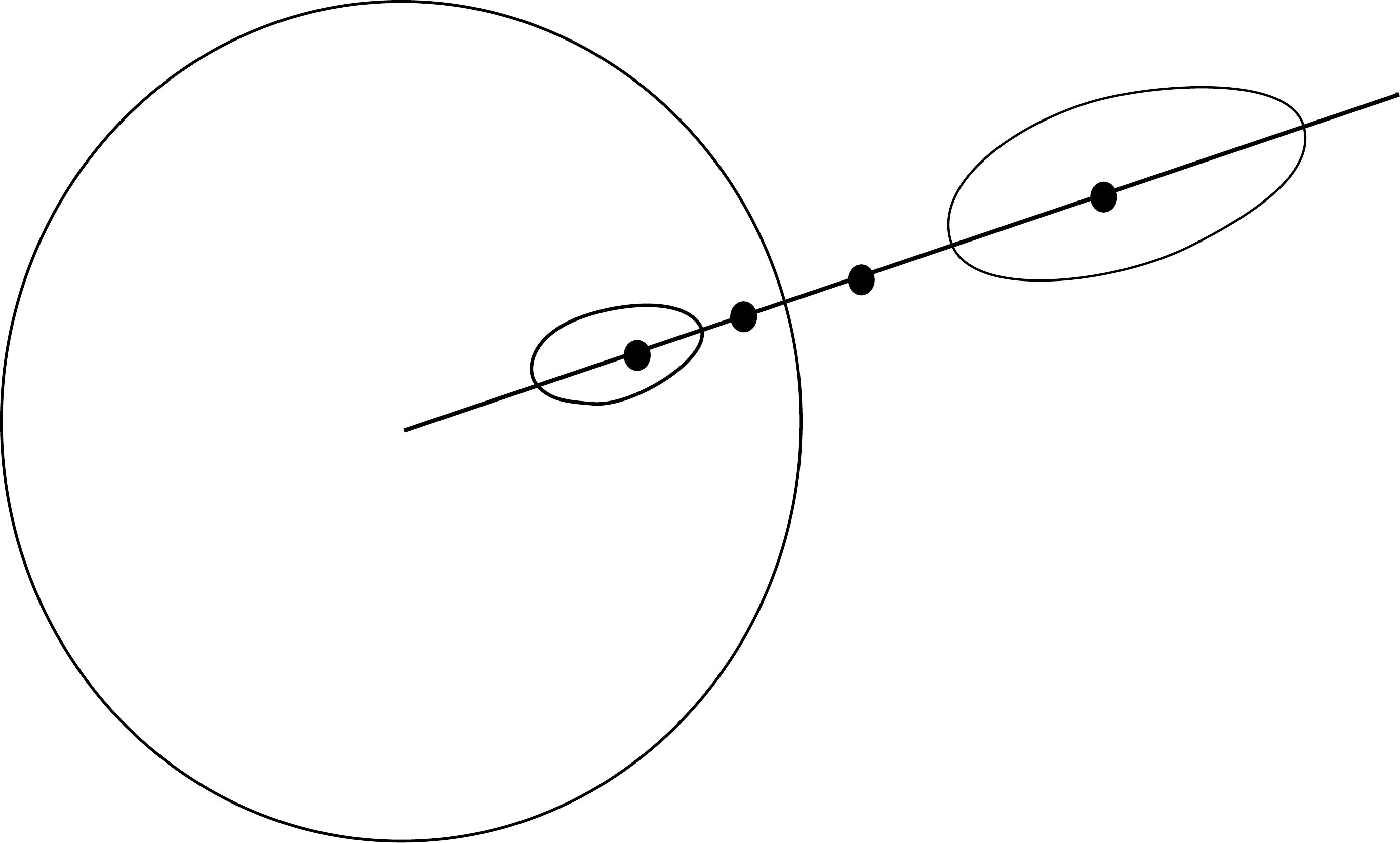}
    \caption{\small{Different configurations of the critical points and the preimages of zero and infinity depending on $|a|$. }} 
    \label{esquemapunts}
\end{figure}

For completeness we describe some features of the dynamics of $B_a$ which depend on the modulus of $a$. The  first thing to consider is whether there is or not a preimage of $\infty$  in $\dis$. This family has a unique pole at $z_{\infty}=1/\overline{a}$ and a unique zero $z_0=a$. Their position, together with the positions of $c_{\pm}$, influence the possible dynamics of $B_a$. We proceed to describe the situation depending on $|a|$ (see Figure~\ref{esquemapunts}).

When $|a|<1$  we have that both critical points $c_{\pm}$ lie on the half ray containing $a$. Moreover, $|c_-|<1$ and $|c_+|>1$. The only pole, $z_{\infty}=1/\overline{a}$ has modulus greater than one. Hence, $B_a:\dis\rightarrow\dis$ is a holomorphic self map of $\dis$ having $z=0$ as a superattracting fixed point. Since, by symmetry, there is no preimage of the unit disk outside the unit circle, $B_a|_{\dis}$ is a degree $4$ branched covering. By  Schwarz Lemma we have that $z=0$ is the only attracting point of $B_a$ in $\dis$ which attracts all points in $\dis$. Summarizing, we have:

\begin{lemma}\label{juliacercle}
If $|a|<1$, $A_a(0)=A_a^*(0)=\dis$ and by symmetry  $A_a(\infty)=A_a^*(\infty)=\widehat{\com}\setminus\overline{\dis}$. Hence, $\mathcal{J}(B_a)=\cercle$.

\end{lemma}

When $|a|=1$ both critical points and the preimages of $0$ and $\infty$ collapse at the point $z=a$, where the function is not formally defined. Everywhere else we have the equality:

$$B_a(z)=z^3\frac{z-a}{1-z/a}=-az^3.$$

\begin{figure}[hbt!]
    \centering
    \subfigure[\scriptsize{Dynamical plane of $B_{3/2i}$}  ]{
    \includegraphics[width=190pt]{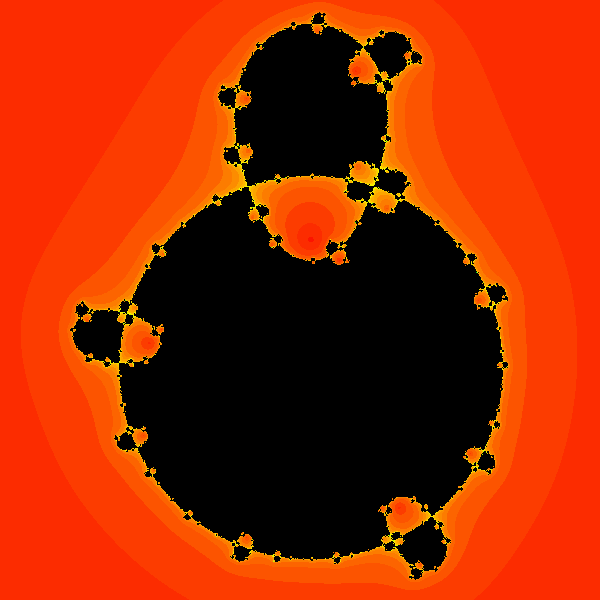}}
    \hspace{0.1in}
    \subfigure[\scriptsize{Dynamical plane of $B_{3/2}$}  ]{
    \includegraphics[width=190pt]{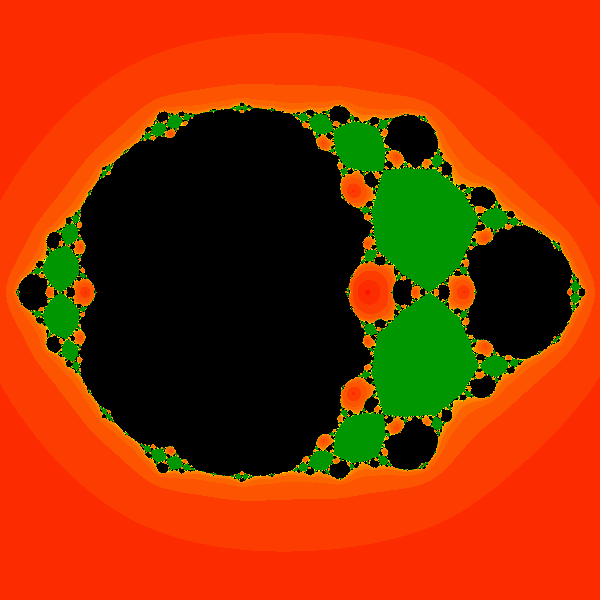}}
    \hspace{0.1in}  
 \subfigure[\scriptsize{Dynamical plane of $B_{1.07398+0.5579i}$}  ]{
     \includegraphics[width=190pt]{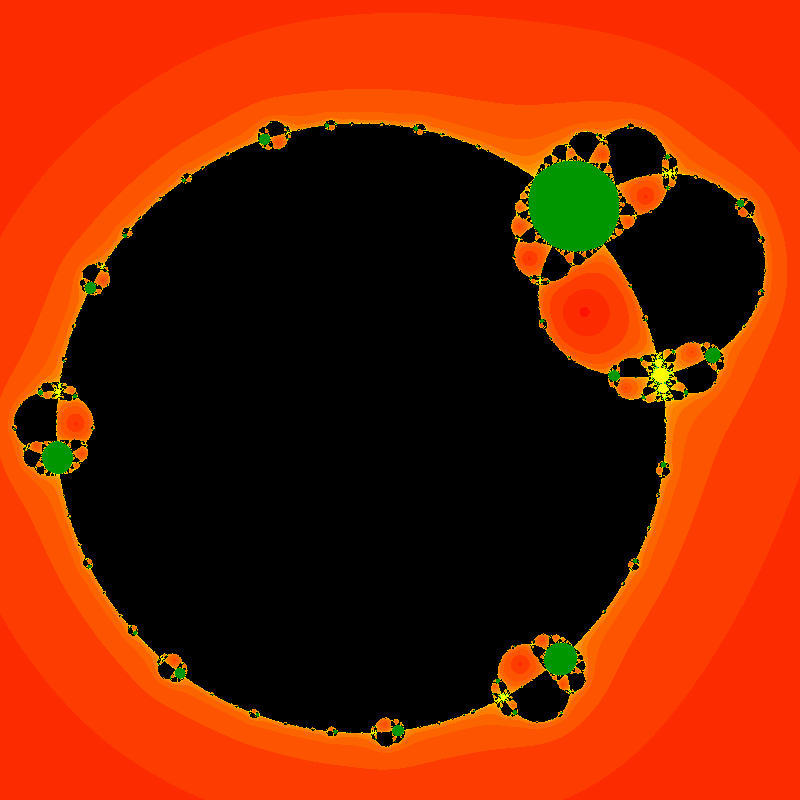}}
    \hspace{0.1in}  
   \subfigure[\scriptsize{Zoom in (c)}  ]{
    \includegraphics[width=190pt]{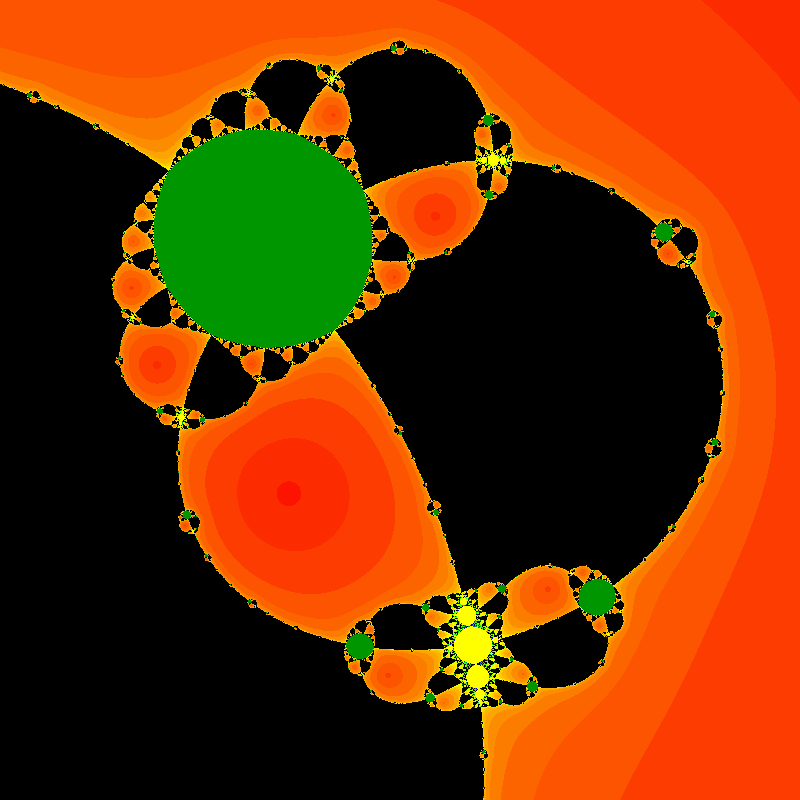}}
    \caption{\small{Dynamical planes of three Blaschke products $B_a$ with $1<|a|<2$. The colours work as follows: a scaling of red if the orbit tends to infinity, black if it tends to zero, green if the orbit accumulates on the cycle $<z_0>$ such that $c_+\in A^*(<z_0>)$ and yellow if there exists a cycle $<w_0>\neq<z_0>$ such that $c_-\in A^*(<w_0>)$ and the orbit accumulates on it. In case (a) there are no other Fatou components than the basins of zero and infinity. In case (b) both free critical points are attracted to a period 2 cycle. In Figures (c) and (d) the critical points are attracted to two different cycles of period 1 (green) and period 4 (yellow), respectively.}} 
    \label{1<a<2}
\end{figure}

For $1<|a|<2$, it follows from (\ref{criticalpoints}) that $c_+=a\cdot k$ and $c_-=a\cdot \bar{k}$ with  $k\in\com$. Hence, by symmetry, $|c_{\pm}|=1$ (see Figure~\ref{esquemapunts}~(b)). The critical orbits lie in $\cercle$ and are not related to each other. The circle map $B_{a}|_{\cercle}$ has no degree defined. Indeed, it can be proven that some points in $\cercle$ have $2$ preimages under  $B_{a}|_{\cercle}$ whereas other points have $4$ preimages.  In Figure \ref{1<a<2}  we show the dynamical planes of three maps $B_a$ with $1<|a|<2$.

\begin{figure}[hbt!]
    \centering
    \subfigure{
    \includegraphics[width=190pt]{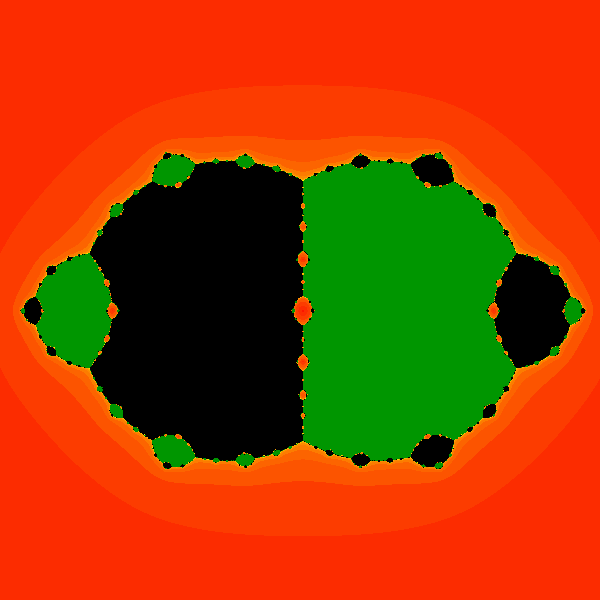}}
    \hspace{0.1in}
    \subfigure{
    \includegraphics[width=190pt]{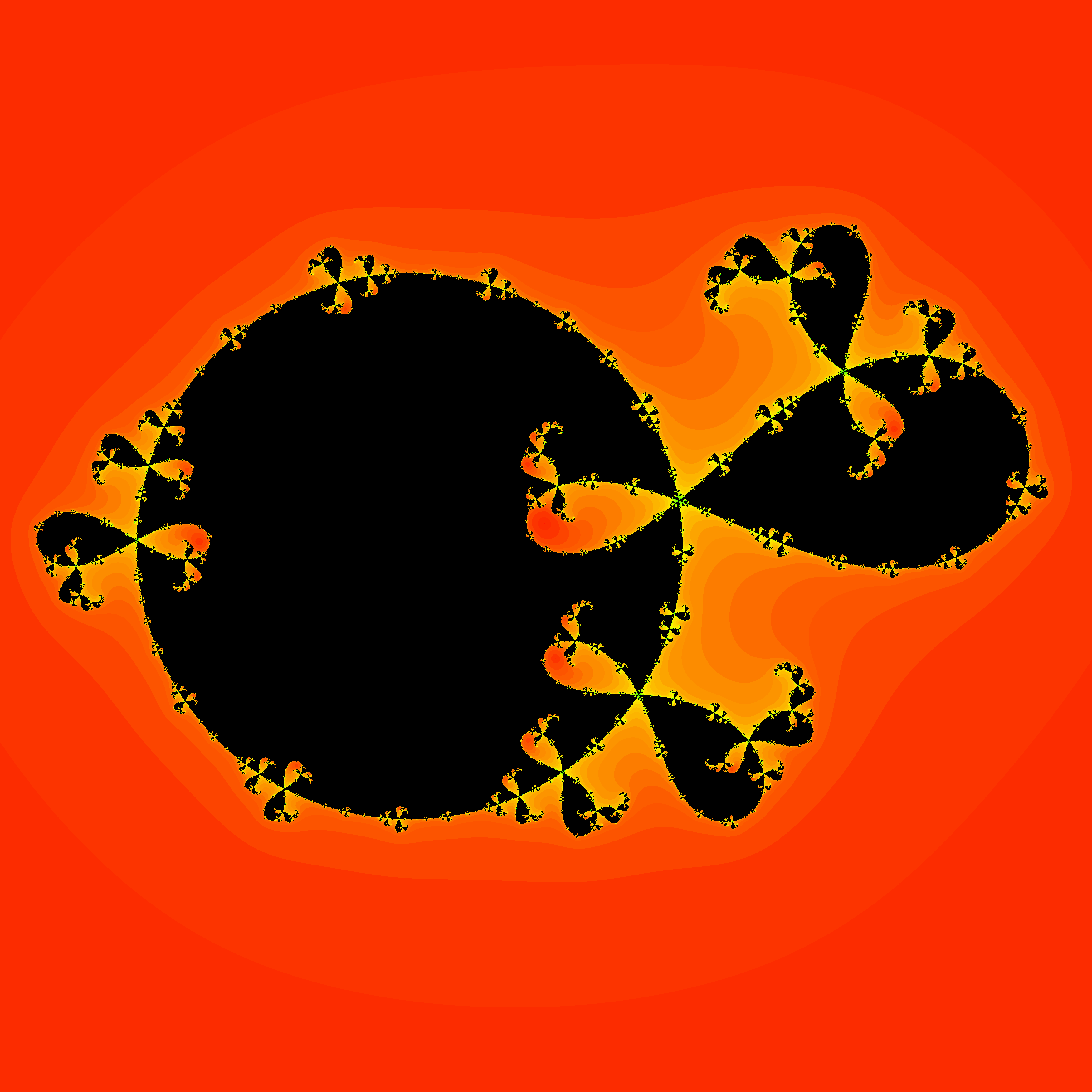}}
  \caption{\small{Dynamical planes of $B_2$ (left) and $B_{a_0}$, where $a_0=1.971917+0.333982i$, (right) . In the left case the point $z=1$ is a superattracting fixed point. The parameter $a_0$ has been chosen numerically so that $B_{a_0}|_{\cercle}$ has no attracting cycle. Colours are as in Figure \ref{1<a<2}. }}   
  \label{a=2}
\end{figure}

When $|a|=2$ we have a unique critical point $c=a/2$ of multiplicity 2 in the unit circle. There are two preimages of $\cercle$ which meet at $c$ (see Figure \ref{esquemapunts} (c)). There may or may not be an attracting or parabolic cycle in $\cercle$ when $|a|=2$. The parameter might be, for example, of Misiurewicz type (i.e.\ the free critical point is preperiodic). In this situation the only Fatou components of $B_a$ are the basins of $z=0$ and $z=\infty$. We also remark for further use that the map $B_a|_{\cercle}$ is 2-to-1. In Figure \ref{a=2} we show the dynamical planes of two maps $B_a$ with $|a|=2$.

\begin{figure}[hbt!]
    \centering
    \subfigure{
    \includegraphics[width=190pt]{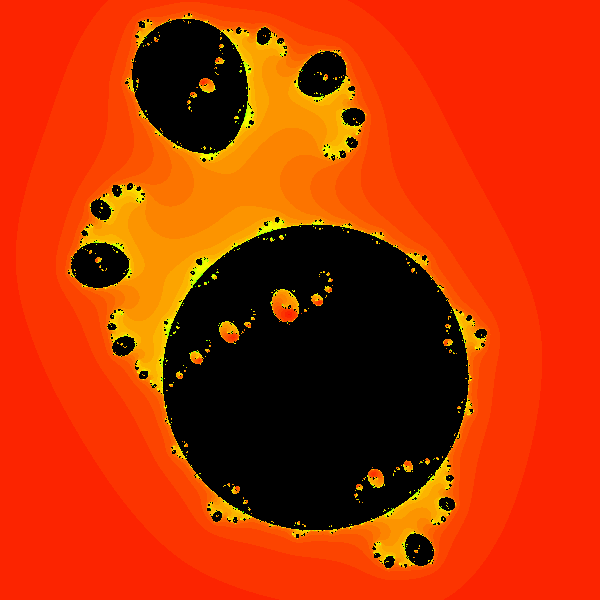}}
    \hspace{0.1in}
    \subfigure{
    \includegraphics[width=190pt]{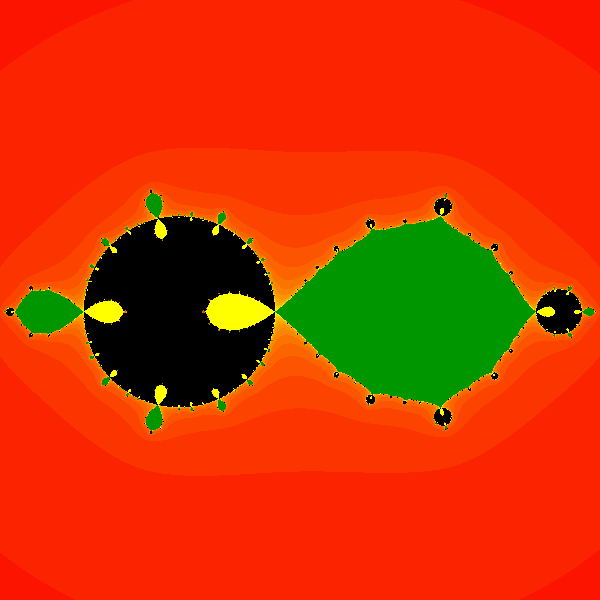}}
  \caption{\small{Dynamical planes of $B_{a_0}$, where $a_0=-0.87+2.05333i$, (left) and $B_4$ (right). In the left case  the critical point $c_+$ belongs to $A^*(\infty)$ and the Julia set is disconnected.  In the right case each free critical orbit accumulates on a different basin of attraction. Colours are as in Figure \ref{1<a<2}.}}   
  \label{a>2}
\end{figure}

When $|a|>2$, as is the case when $|a|<1$,  we have that both  critical points $c_{\pm}$ lie on the half ray containing $a$ and are symmetric with respect to $\cercle$. In this case we have two disjoint preimages of the unit circle: one of them inside $\dis$, surrounding the pole $z_{\infty}$, and the symmetric one outside surrounding the zero $z_0=a$ (see Figure \ref{esquemapunts} (d)). As in the case $|a|=2$, $B_a|_{\cercle}$ is 2-to-1.   
In Figures \ref{a>2} and \ref{hyperbolicpic} (a) we show the dynamical planes of three maps $B_a$ with $|a|>2$.

\subsection{Connectivity of the Julia set: proof of Theorem A}\label{sconn}

The goal of this section is to prove Theorem A. Notice that statement (a) has already been proven in  Lemma \ref{juliacercle}.  In Proposition \ref{conninfty} we prove statement (b).

\begin{propo}\label{conninfty}
Let  $B_a$ be as in (\ref{blasformula}) and suppose $|a|>1$. Then, the connected components of $A(\infty)$ and $A(0)$ are simply connected if and only if $c_+\notin A^*(\infty)$. 
\end{propo}
\proof
 By symmetry, the connected components of $A(0)$ are simply connected if and only if the ones of $A(\infty)$ are. Therefore, we focus on the simple connectivity of $A(\infty)$. By means of the Riemann-Hurwitz formula (Theorem \ref{riemannhurwitz}) and invariance of $\cercle$, the connected components of $A(\infty)\setminus A^*(\infty)$ are simply connected if and only if $A^*(\infty)$ is simply connected since any connected component of $A(\infty)\setminus A^*(\infty)$ can have at most one critical point (see Corollary \ref{connpreim}). Therefore, it is sufficient to prove that $A^*(\infty)$ is simply connected if and only if $c_+\notin A^*(\infty)$. Let us consider the B\"ottcher coordinate of the superattracting fixed point $z=\infty$ (see \cite{Mi1}, Theorem 9.1). If there is no extra critical point in $A^*(\infty)$, the B\"ottcher coordinate can be extended until it reaches $\partial A^*(\infty)$ and $A^*(\infty)$ is simply connected. If it does contain an extra critical point, the B\"ottcher coordinate can only be extended until it reaches the critical point  (see \cite{Mi1}, Theorem 9.3). Let $U$ be the maximal domain of definition of the B\"ottcher coordinate at $\infty$. Then, either $\partial U$ consists of the union of two topological circles, say $\gamma_{\pm}$ , which are joined in a unique point which is the critical point, or $\partial U$ is a topological circle containing the critical point. If it is the last case, there is an extra preimage of $B_a(U)$ attached to the critical point. Hence, $ A^*(\infty)$ would be mapped 4 to 1 onto itself. This is not possible since $B_a$ is of degree 4 and the only pole $z_{\infty}$ is inside the unit disk and hence does not belong to $A^*(\infty)$. Let $V_+$ and $V_-$ be the disjoint simply connected regions bounded by $\gamma_+$ and $\gamma_-$. The result follows by noticing that both $B_a(V_+)$ and $B_a(V_-)$ contain $\widehat{\com}\setminus B_a(U)$ and, hence, both of them  contain Julia set since $\mathcal{J}(B_a)$ is not empty and $U$ is contained in $\mathcal{F}(B_a)$.
\endproof

We now begin the proof of statement (c). In propositions \ref{noHR} and \ref{connother} we prove that all periodic Fatou components other than $A^*(0)$ or $A^*(\infty)$ are simply connected.

\begin{propo}\label{noHR}
Let  $B_a$ be as in (\ref{blasformula}). Then  $B_a$ has no Herman Rings.
\end{propo}
\proof
M.~Shishikura \cite{Shi} proved that if a rational map has a Herman ring, then it has two different critical points whose orbits accumulate on the two different components of the boundary of it. If $|a|\leq 1$, we have that the Julia set satisfies $\mathcal{J}(B_a)=\cercle$ (see Lemma \ref{juliacercle}), so $B_a$ cannot have Herman rings. If $1<|a|\leq2$, the two critical orbits lie in $\cercle$ and, hence, there can be no Herman rings.

  We focus now on the case $|a|>2$.  By Shishikura's result and symmetry, a cycle of Herman rings would have components both inside and outside the unit disk. Hence, it would have at least one component in the preimage  of the unit disk $\Omega_e=B_a^{-1}(\dis)\setminus\dis$ and another one in the preimage of the complementary of the unit disk $\Omega_i=B_a^{-1}(\com\setminus\overline{\dis})$ (see Figure~\ref{esquemapunts}~(d)). Recall that $\Omega_e$ is a simply connected set disjoint from $\cercle$. Moreover, all its preimages are bounded, none of them can intersect the unit circle and all of them are simply connected by Corollary \ref{connpreim}. Every component of the cycle of Herman rings is contained in a preimage of $\Omega_e$ of some order $n\geq 0$. We claim that such a cycle must have a component which surrounds either the unit disk or $\Omega_e$. If this is so, this component cannot be contained in a simply connected preimage of $\Omega_e$, which leads to a contradiction.

 Let $\mathcal{I}(z)=1/\overline{z}$ be the reflection with respect to $\cercle$. To prove the claim observe that, due to symmetry, if $A$ is a component of the cycle of Herman rings, then so  is $\mathcal{I}(A)$. Moreover, since infinity is a superattracting fixed point, all  components are bounded and at least one of them, say $A'$, surrounds the pole $z_{\infty}$ (by the Maximum Modulus Principle).  Recall that $z_{\infty}$ is contained in $\Omega_i$  and that, again by symmetry, $\mathcal{I}(\Omega_i)=\Omega_e$. Then, either  $A'$ surrounds the unit disk or surrounds $\Omega_i$ or is contained in $\Omega_i$. In the first case we are done. In the second case $\mathcal{I}(A')$ surrounds $\Omega_e$ and we are also done. In the third case, $B_a(A')$ separates infinity and the unit disk and, hence, surrounds the unit disk. This finishes the proof.

\endproof

\begin{propo}\label{connother}
Let  $B_a$ be as in (\ref{blasformula}) with $|a|\geq 2$. Let $<z_0>$ be an attracting, superattracting or parabolic $p$-cycle of $B_a$ other than $\{0\}$ or $\{\infty\}$. Then $A^*(<z_0>)$ is simply connected.
\end{propo}

\proof

\underline{Case 1:} First we consider the case in which each connected component of the immediate basin of attraction contains at most one critical point (counted without multiplicity). For the attracting case consider a linearizing  domain $\mathcal{A}$ coming from K{\oe}nigs linearization around $z_0$ (see \cite{Mi1}, Theorem 8.2). The subsequent preimages $U_n$ defined as the components of $B_a^{-n}(\mathcal{A})$ such that $z_{-n}\in U_n$, contain at most one critical point and are hence simply connected by Corollary \ref{connpreim}. The result follows since the nested subsequence of preimages $\{U_{np}\}$ covers $A^*(z_0)$. The parabolic case follows similarly taking a petal instead of a linearizing domain (see \cite{Mi1}, Theorem 10.7) whereas in the superattracting case we may use a B\"otcher domain (see \cite{Mi1}, Theorem 9.1).

\underline{Case 2:} Now consider the case in which one connected component, say $A^*(z_0)$, of the immediate basin of attraction contains the two free critical points. This excludes the case $|a|=2$ (see Section \ref{introblas}). Without loss of generality we assume that $z_0$ is a fixed point. Indeed, the first return map from  $A^*(z_0)$ onto itself has no other critical points since the other components of the immediate basin of attraction contain none.

Due to symmetry of the critical orbits, the fixed point $z_0$ lies in $\cercle$. Hence, $A^*(z_0)$ intersects $\cercle$, which is invariant. If $z_0$ is attracting, take the maximal domain $\mathcal{A}$ of the K{\oe}nigs linearization (see \cite{Mi1}, Lemma 8.5). Its boundary $\partial \mathcal{A}$ contains, due to symmetry, the two critical points. Each critical point has a different simply connected preimage of $B_a(\mathcal{A})$ attached to it. Now consider $V=B_a^{-1}(\mathcal{A})$. The map $B_{a|V}:V\rightarrow \mathcal{A}$ is proper and of degree 3 since $z_0$ has three different preimages. Given that $V$ contains exactly 2 critical points and $B_{a|V}$ is of degree 3, it follows from the Riemann-Hurwitz formula (see Theorem \ref{riemannhurwitz}) that $V$ is simply connected. Using the same reasoning all of its preimages are simply connected. Finally, since $A^*(z_0)$ is covered by the nested sequence of simply connected preimages of $\mathcal{A}$, we conclude that $A^*(z_0)$ is simply connected. The parabolic case is done similarly taking $\mathcal{P}$ to be the maximal invariant petal (see \cite{Mi1}, Theorems 10.9 and 10.15). Notice that, due to symmetry, for $|a|>2$ we cannot have a superattracting cycle of local degree  2 with an extra critical point in $A^*(z_0)$ .

\endproof

 We now finish the proof of statement (c). Assume that there exists a periodic Fatou component other than $A^*(0)$ and $A^*(\infty)$. Then, such periodic Fatou component has a critical point related to it. Indeed, if it is Siegel  disk, there is critical point whose orbit accumulates on its boundary (see \cite{Shi}). If it is the basin of attraction of an attracting, superattracting or parabolic cycle $<z_0>$, there is  a critical point in $c\in A^*(<z_0>)$  (cf.\ \cite{Mi1}). Therefore, there is at most one unoccupied critical point. Hence, by means of the Corollary \ref{connpreim} of the Riemann-Hurwitz formula, any preperiodic Fatou component eventually mapped to a periodic component other that $A^*(\infty)$ or $A^*(0)$ is also simply connected.
 
 The final statement of the theorem holds since the Julia set of a rational map is connected if and only if all connected Fatou components are simply connected.

\section{Parameter plane of the Blaschke family}\label{parameterplane}

\begin{figure}[hbt!]
\centering
\includegraphics[width= 10cm]{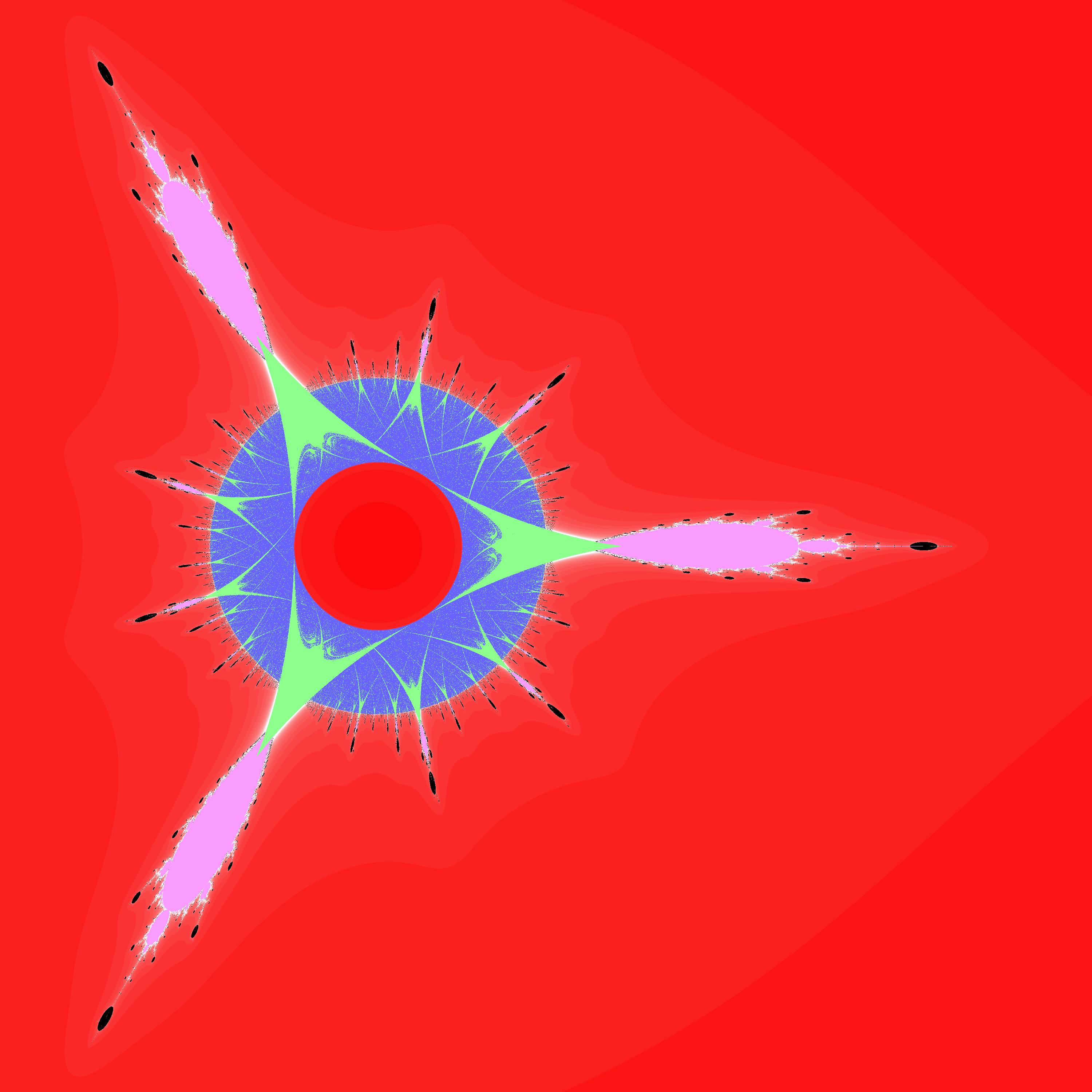}
\caption{\small Parameter plane of the Blaschke family $B_a$. The colours are as follows: red if $c_+\in A(\infty)$, black if  $c_+\in A(0)$, green if $O^+(c_+)$ accumulates on a periodic orbit in $\cercle$, pink if $O^+(c_+)$ accumulates in a periodic orbit not in $\cercle$ and blue in any other case. The inner red disk corresponds to the unit disk. }
\label{paramblash}
\end{figure}

The aim of this section is to study the parameter plane of the Blaschke family $B_a$. Figure \ref{paramblash} shows the result of iterating the critical point $c_+$. Since the two critical orbits of $B_a$ are related by symmetry unless $1<|a|<2$, this information suffices also for $c_-$ everywhere else. Indeed, if $1<|a|<2$, the critical orbits may have completely independent behaviour. (see Figure \ref{1<a<2} (c), (d)) 

The next lemma explains the observed symmetries on the parameter plane. Its proof is straightforward.
\begin{lemma}\label{symm}

Let $\xi$ be a third root of the unity. Then, $B_a$ and $B_{\xi a}$ are  conjugate by the conformal map $\tau(z)=\bar{\xi}z$. Moreover, $B_a$ and $B_{\bar{a}}$ are conjugate by the anticonformal map $\widetilde{\mathcal{I}}(z)=\overline{z}$.
\end{lemma}

\subsection{Hyperbolic parameters}\label{hypbeh}

We say that a parameter $a\in\com$ is \textit{hyperbolic} if $B_a$ is hyperbolic, that is, if both free critical points $c_{\pm}=c_{\pm}(a)$ are attracted to attracting cycles. Recall that the two critical orbits are symmetric except for $1<|a|\leq2$, in which case they belong to $\cercle$. Hence, if one fixed point is attracted to the superattracting fixed point $z=0$ (resp.\ $z=\infty$) the other one is attracted to $z=\infty$ (resp.\ $z=0$). Parameters for which this happens are called \textit{escaping parameters}. We denote the set of escaping parameters by $\mathcal{E}$. We shall denote the set of non-escaping parameters by $\mathcal{B}$.
Observe that $\mathcal{B}$ is not the connectedness locus of $B_a$. Indeed, in view of Theorem A, parameters in $\mathcal{E}$ have a connected Julia set if the critical points belong to $A(0)\setminus A^*(0)$ or $A(\infty)\setminus A^*(\infty)$. On the other hand, some Julia sets for $1<|a|<2$ (which must belong to $\mathcal{B}$) may, a priori, be disconnected (although we believe that this is never the case). Notice also that, if $|a|=1$, the products $B_a$ degenerate to degree $3$ polynomials without free critical points in which case $a$ is neither escaping nor non-escaping.

\begin{lemma}
If $|a|<1$ then $a\in\mathcal{E}$. If $1<|a|\leq 2$ then $a\in\mathcal{B}$. The non-escaping set $\mathcal{B}$ is bounded.
\end{lemma}
\proof
The first two statements are already proven.  To prove the third one we have to see that, if $|a|$ is big enough, then the parameter $a$ is escaping. First we prove that, if $|z|>\lambda(|a|+1)$ with $\lambda\geq 1$, then $|B_a(z)|>\lambda|z|$. It follows from the previous hypothesis that $|z-a|>\lambda$ and that $|z|^2>|z|(|a|+1)>|1-\overline{a}z|$. Therefore, we have
$$|B_a(z)|=|z|^3\frac{|z-a|}{|1-\overline{a}z|}>|z|^3\frac{\lambda}{|z|^2}=\lambda|z|.$$

To finish the proof notice that, as $|a|$ tends to infinity, the critical point $c_+(a)$ tends to $a\cdot 2/3$. Consequently, it can be checked that the modulus of the critical value $v_+=B_a(c_+(a))$ grows as $\mathcal{C} |a|^2$ for some $\mathcal{C}>0$ and, for $|a|$ big enough, $|v_+|>\lambda(|a|+1)$ with $\lambda>1$. We conclude that $|B_a^n(v_+)|\rightarrow\infty$ when $n\rightarrow\infty$. Therefore, for $|a|$ big enough, $a\in\mathcal{E}$.

\endproof

\begin{figure}[hbt!]
    \centering
    \subfigure{
    \def\svgwidth{190pt}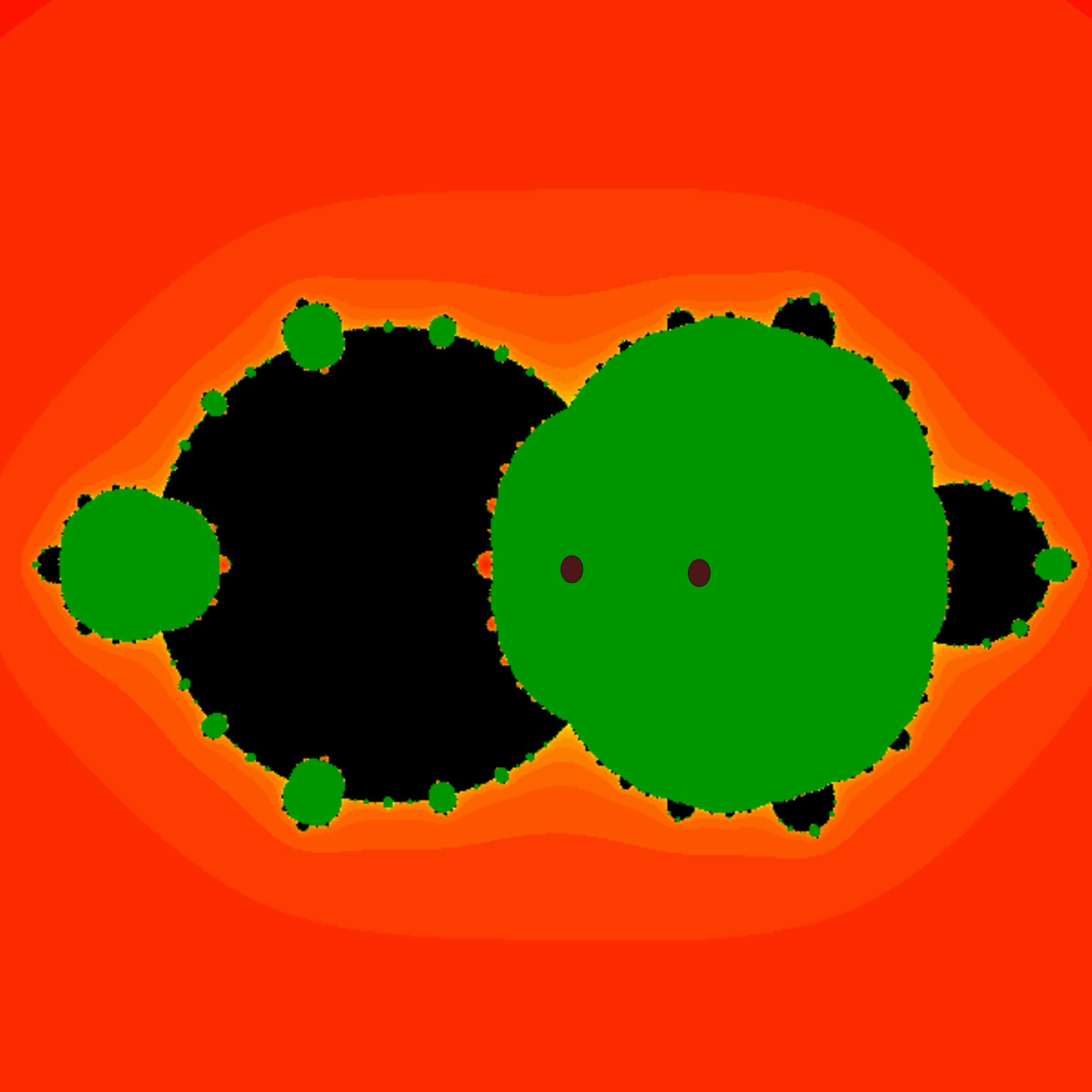}
    \hspace{0.1in}
 \subfigure{
     \def\svgwidth{190pt}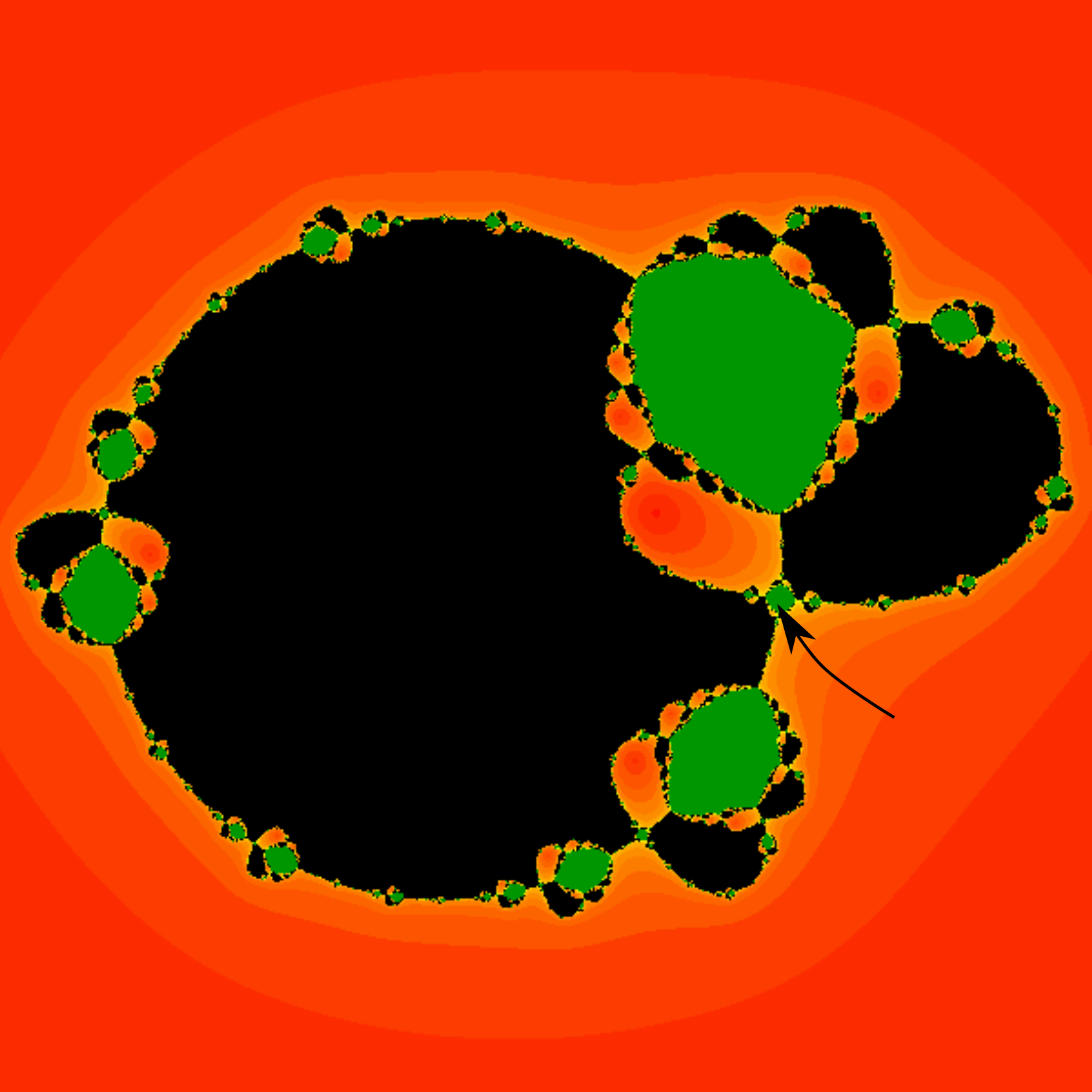}
    \hspace{0.1in}  

    \caption{\small{Dynamical planes of the Blaschke products $B_{2.5}$ (left) and $B_{1.52+0.325i}$ (right). The colours are as in Figure \ref{1<a<2}. The left case corresponds to an adjacent parameter $a$ in a tongue (see Definition \ref{deftongue}).  The right case corresponds to a capture parameter. }} 
    \label{hyperbolicpic}
\end{figure}

Following  \cite{Re} and  \cite{Mi3}, we classify hyperbolic non-escaping parameters as follows.
\begin{defin}
Let $B_a$ be a hyperbolic Blaschke product such that  $a\in \mathcal{B}$, $|a|\neq 2$. We say that the map is:
\begin{enumerate}[(a)]
\item \textit{Adjacent} if the critical points belong to the same component of the immediate basin of attraction of an attracting cycle (see Figure \ref{hyperbolicpic} (left)).
\item \textit{Bitransitive} if the  critical points belong to different components of the same immediate basin of attraction of an attracting cycle (see Figure \ref{1<a<2} (b)).
\item \textit{Capture} if one of the critical points belongs to the immediate basin of attraction of an attracting cycle and the other one belongs to a preperiodic preimage of it (see Figure \ref{hyperbolicpic} (right)).
\item  \textit{Disjoint} if the  critical points belong to the immediate basin of attraction of two different attracting cycles (see Figure \ref{1<a<2} (c) and (d), Figure \ref{a>2} (right) and Figure \ref{bitransdisjointswap} (right)).
\end{enumerate}
\end{defin}

We have omitted the parameters $a$ with $|a|=2$ in the previous definition since, for them, the two free critical points collide in a single one.

Because of the non holomorphic dependence in the parameter $a$, this parameter slice presents some special hyperbolic components which we call \textit{tongues} in analogy to the Arnold tongues which appear in the Arnold standard family of perturbations of rigid rotations (see \cite{Ar}, \cite{Fa2} and \cite{Ge}). In this settings they cannot be defined as the level sets of the rotation number (since there is no rotation number to speak about), but instead, we may use the equivalent characterization of the existence of an attracting periodic orbit in $\cercle$ with given combinatorics. Tongues for the double standard map were studied previously by M.~Misiurewicz and A.~Rodrigues \cite{MiRo1, MiRo2} and A.~Dezotti \cite{De}. They are defined as follows.

\begin{defin}\label{deftongue}
 We define $\mathcal{T}$ as the set of parameters $a$, $|a|\geq 2$, such that $B_a|_{\cercle}$ has an attracting cycle. The connected components of $\mathcal{T}$ are called \textit{tongues}.
\end{defin} 

Tongues can be be classified according to the combinatorics of the attracting cycle, which can be described via the semiconjugacy between $B_a|_{\cercle}$ and the doubling map (see Lemma \ref{semiconj}). A detailed study of $\mathcal{T}$ will appear in a forthcoming paper.

 We proceed now to describe in which regions the different types of non-escaping hyperbolic parameters belong.

\begin{teor}
Let $a\in \mathcal{B}$ with $|a|\neq 2$. Then,

\begin{enumerate}[(a)]\label{clashyp}
\item If $a$ is an adjacent parameter, either $1<|a|<2$ or it belongs to a  tongue. Conversely, any parameter $a$ belonging to a tongue  is adjacent. 
\item If $a$ is a bitransitive parameter, then either $1<|a|<2$ or $c_+$ and $c_-$ enter and exit the unit disk infinitely many times.
\item If $a$ is a capture parameter, then $1<|a|<2$.
\item If $a$ is a disjoint parameter, then either $1<|a|<2$ or $|a|>2$. In the later case the orbits of the two attracting cycles are symmetric with respect to the unit disk and, hence, have the same period. Moreover,   if the multiplier of one attracting cycle is $\lambda$, the multiplier of the other attracting cycle is $\overline{\lambda}$. 
\end{enumerate}

\end{teor}

\proof
We begin with statement (a).  If $a$ is an adjacent parameter, then both critical points belong to the same component of the immediate basin of attraction of a periodic cycle. Then, either $1<|a|<2$ or $|a|>2$ and both critical points are attracted to an attracting cycle in $\cercle$. In this last case, $a$ belongs, by definition, to a tongue. The converse holds by symmetry.

To prove statement (b) notice that, in the bitransitive case, the immediate basin of attraction of the attracting cycle  which the critical points are attracted to has at least two different connected components. If $|a|>2$, by symmetry, at least one is contained in the unit disk and another one is on its complement. Thus, the critical orbits enter and exit $\dis$ infinitely many times.

Statement (c) follows directly from the fact that, for $|a|\geq2$, the critical orbits are symmetric.

The first part of (d) follows from symmetry. In order to see that the attracting cycles have conjugated multipliers, we conjugate $B_a$ via a M\"obius transformation $M$ to a rational map $\widetilde{B}_a$ that fixes the real line. The result follows then from the fact that $\widetilde{B}_a'(\bar{z})=\overline{\widetilde{B}'_a(z)}$ and that $M$ preserves the multiplier of the periodic cycles.

\endproof

\subsection{Relation with cubic polynomials}\label{surgery}

\begin{figure}[hbt!]
    \centering
    \subfigure{
    \includegraphics[width=190pt]{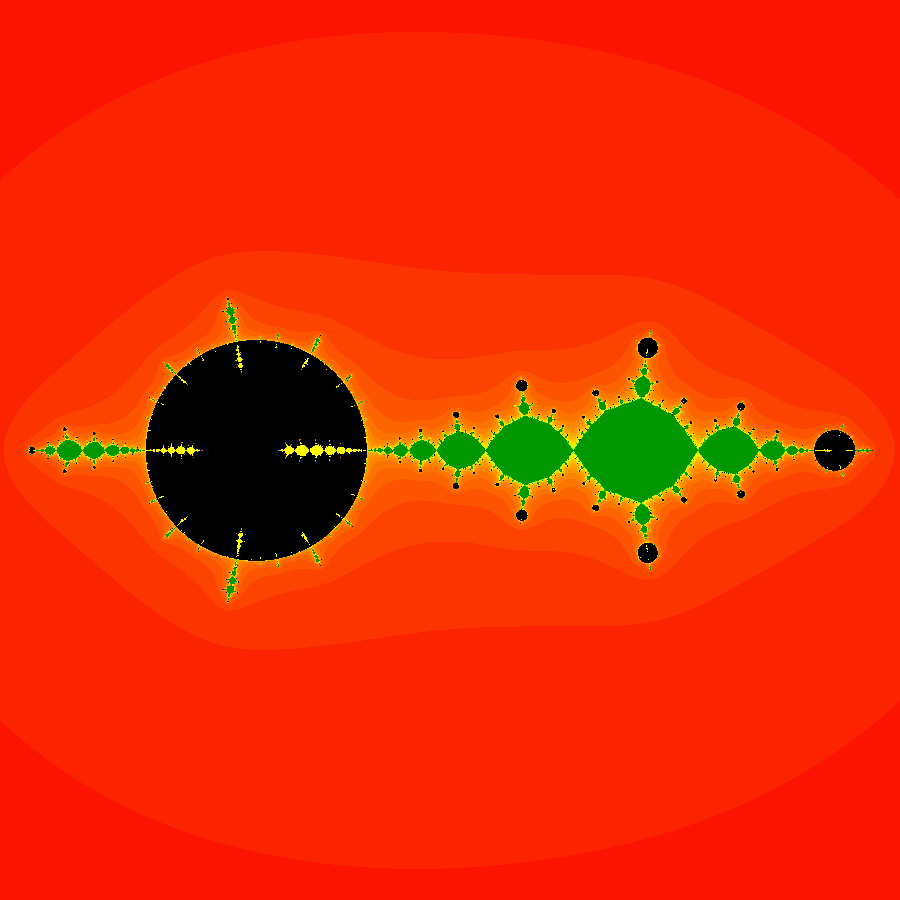}}
    \hspace{0.1in}
    \subfigure{
    \includegraphics[width=190pt]{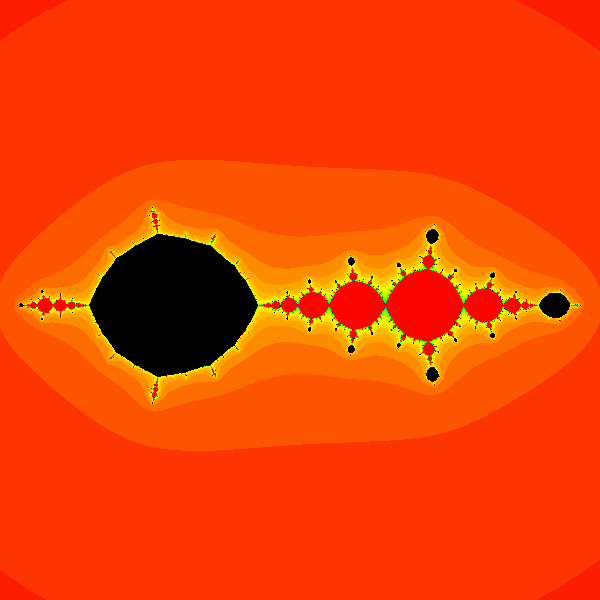}}
    \caption{\small{Dynamical planes of the Blaschke product $B_{5.25}$ (left) and the cubic polynomial $M_{-5.5}$ (right). The black regions of both figures correspond to the basins of attraction of the superattracting fixed points $z=0$. The cubic polynomial shows, in red, the basin of attraction of a period two attracting cycle. The Blaschke product has two different attracting cycles of period two ($5.25$ is a disjoint parameter). One outside the unit disk (green) and the other one  inside (yellow). }} 
    \label{blashmilfig}
\end{figure}

\begin{figure}[hbt!]
    \centering
       \includegraphics[width=250pt]{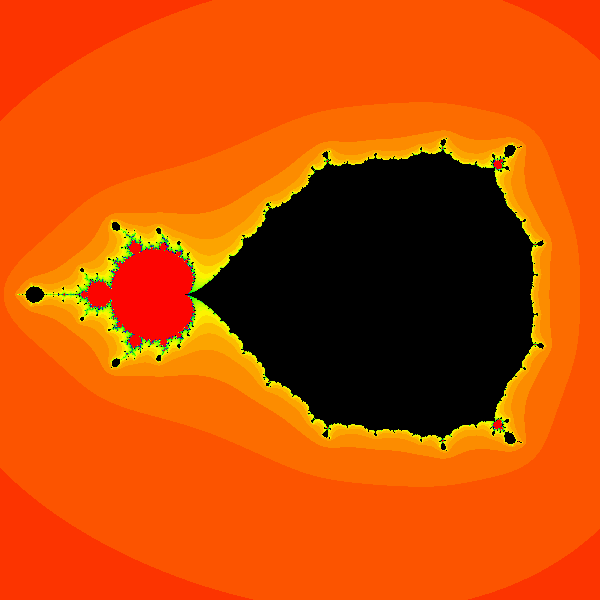}
   
 \caption{\small Parameter plane of the polynomials $M_b$. Colours go as follows: black if the free critical orbit   tends to the superattracting cycle $z=0$ and  red if it tends neither to $z=0$ nor to $z=\infty$. The scaling from green to orange corresponds to parameters for which the critical orbit tends to $z=\infty$.}
 \label{parametercubic}
\end{figure}

In this section we introduce a quasiconformal surgery which relates, under certain conditions, the dynamics of a Blaschke product $B_a$ outside the unit disk with the ones of a cubic polynomial of the form $M_b(z)=bz^2(z-1)$ with $b\in\com$ (see Figure \ref{blashmilfig}). They have $z=0$ as a superattracting fixed point and a second critical point, $c=2/3$, which is free. This one dimensional slice (or a cover thereof) was introduced by Milnor in 1991 in a preliminary version of \cite{Mi2}. They have also been studied by P.~Roesch \cite{Ro}, among others.  In Figure \ref{parametercubic} we show its parameter plane.

We proceed to introduce the quasiconformal surgery which relates the Blaschke products $B_a$ with the cubic polynomials $M_b$ (cf.\ \cite{Pe}). For an introduction to the tools used in quasiconformal surgery we refer to  \cite{Ah} and \cite{BF}.  The idea of the surgery is to ``glue'' the map $R_2(z)=z^2$ inside $\dis$ keeping the dynamics of $B_a$ outside $\dis$ whenever the parameter $a$ is such that $B_{a}|_{\cercle}$ is quasisymmetrically conjugate to the doubling map $\theta\rightarrow 2\theta\; (\modul 1)$. 

More precisely, we  restrict to the set of parameters $a$ such that $|a|\geq 2$. For these parameters, $B_a|_{\cercle}$ is a degree 2 cover of $\cercle$ and is hence semiconjugate to the doubling map  by a non decreasing continuous map $h_a$, not necessarily surjective (see Lemma \ref{semiconj}). We know, by Theorem  \ref{conjpe}, that if $|a|>2$ and the circle map $B_a|_{\cercle}$ has neither attracting nor parabolic cycles, then $h_a$ is a quasisymmetric isomorphism of the circle. Therefore, we define  $\mathcal{X}$ to be the set of parameters $a$, $|a|\geq 2$, such that $h_a$ is a quasisymmetric conjugacy between $B_a|_{\cercle}$ and $R_2$.

Let  $a\in \mathcal{X}$. The map $h:=h_a$ is quasisymmetric and conjugates $B_a|_{\cercle}$ with the doubling map. Since $h$ is quasisymmetric, it extends to a quasiconformal map $H:\overline{\dis}\rightarrow\overline{\dis}$ (see \cite{BeAh}, \cite{DoEa}, cf.\ \cite{BF}). We define the model map as follows:

$$F(z)=\left\{\begin{array}{lcl}
B_a(z) &  \mbox{ for } & |z|>1\\
H^{-1}\circ R_2\circ H(z) &  \mbox{ for } & |z|\leq 1.
\end{array}\right. $$

\begin{propo}\label{surgerymilnor}

Let $a\in\mathcal{X}$. Then, there exists $b\in\com$ and a quasiconformal map $\psi:\com\rightarrow\com$ such that $\psi\circ F\circ \psi^{-1}=M_{b}$, where $M_b(z)=b^2z(z-1)$.

\end{propo}
\proof

The map $F$ is quasiregular since it is continuous in $\widehat{\com}$, holomorphic outside $\overline{\dis}$ and locally quasiconformal in $\dis\setminus\{0\}$. Its topological degree is 3 since gluing the map $z\rightarrow z^2$ in $\overline{\dis}$ decreases the degree of $B_a$ in 1. Indeed, recall from Section \ref{introblas} that $B_a$ has three preimages of $\dis$ and one preimage of $\com\setminus\dis$ in $\dis$. Instead, $F$ has only two preimages of $\dis$ and none of $\com\setminus\dis$ in $\dis$.

We now define an $F$-invariant almost complex structure $\sigma$, i.e, an almost complex structure such that $F^*\sigma=\sigma$, as

$$\sigma=\left\{\begin{array}{lcl}
H^*\sigma_0 &  \mbox{ on } & \mathbb{D}\\
(F^m)^*(H^*\sigma_0) &\mbox{ on } &  F^{-m}(\mathbb{D})\setminus F^{-m+1}(\mathbb{D}),\mbox{ for } m\geq1\\
\sigma_0 &   \mbox{otherwise,} &
\end{array}\right. $$

\noindent where $\sigma_0$ denotes the standard complex structure and $^*$ denotes the pullback operation. By construction, $\sigma$ has bounded dilatation. Indeed, $\sigma$ is, in $\dis$, the pull back of $\sigma_0$ by a quasiconformal map. Everywhere else either we pull back $\sigma|_{\dis}$ by a holomorphic map (so we do not increase the dilatation) or we use the standard complex structure.

Let $\psi$ be the integrating map of $\sigma$ given by the Measurable Riemann Mapping Theorem (see \cite{Ah} p.\ 57, \cite{BF} Theorem 1.28) such that $\psi(H^{-1}(0))=0$, $\psi(\infty)=\infty$ and $\psi(c_+)=2/3$. Then, the following  diagram commutes.

$$
\begin{CD}
( \com, \sigma) @>F>> ( \com, \sigma) \\
@VV\psi V @VV\psi V\\
(\com, \sigma_0) @>\psi\circ F \circ \psi^{-1}>> (\com, \sigma_0).
\end{CD}
$$

The composition $\psi\circ F \circ \psi^{-1}$ is a quasiconformal map preserving the standard complex structure and, by Weyl's Lemma (see \cite{Ah} p.\ 16, \cite{BF} Theorem 1.14), $\psi\circ F \circ \psi^{-1}$ is  holomorphic. Since this map has topological degree 3 and no poles it is a cubic polynomial. By the chosen normalization, $z=0$ is a superattracting fixed point and $z=2/3$ is a critical point, hence $\psi\circ F \circ \psi^{-1}=M_{b}$ for some $b\in\com$.
\endproof

Since $F=B_a$ outside $\dis$, it follows that $\psi$ is a quasiconformal conjugacy between $M_{b(a)}$ and $B_a$ on this region, and so are all iterates on orbits which never enter $\dis$. Therefore, for all parameters $a\in\mathcal{X}$ such that the orbit of the exterior critical point $\mathcal{O}(c_+)$ never enters $\overline{\dis}$ all relevant dynamics are preserved (see Figure \ref{blashmilfig}). 
The study of the parameters $a\in\mathcal{X}$ such that  $\mathcal{O}(c_+)$ meets $\overline{\dis}$ is done in the next section. 

The surgery described above above defines a map $\Phi:\mathcal{X}\rightarrow\mathcal{Y}$ between the subset $\mathcal{X}$ of the parameter plane of $B_a$ (see Figure~\ref{paramblash}) and a subset $\mathcal{Y}$ of the parameter plane of $M_b$ (see Figure~\ref{parametercubic}). By Theorem \ref{conjpe} and the fact that every parameter $a$ with $|a|>2$ such that $B_a|_{\cercle}$ has a parabolic cycle belongs to the boundary of a tongue (which is to be proven in a forthcoming paper), the set $\mathcal{X}$ includes all parameters with $|a|>2$ such that $a$ does not belong to the closure of any tongue. The image set $\mathcal{Y}$ does not include any parameter $b$ in the main capture component (i.e.\ the set of parameters for which the basin of $0$ contains both critical orbits). We conjecture that $\Phi$ is a degree 3 cover between $\mathcal{X}$ and $\mathcal{Y}$. Another application of this surgery construction will be explained in Section \ref{parametrizehyp} (see Proposition \ref{jordancurveattr}).

\subsection{Swapping regions: Proof of Theorem B}\label{polantipol}

\begin{defin}
We say that a parameter $a$, $|a|>2$, is a \textit{swapping parameter} if the exterior critical point $c_+$ eventually falls under iteration in $\dis$ (or equivalently if $c_-$ eventually falls in $\widehat{\com}\setminus\overline{\dis}$). 
A maximal open connected set of swapping parameters is called a \textit{swapping region}.
\end{defin}

\begin{figure}[hbt!]
    \centering
    \subfigure[\scriptsize{The Tricorn.} ]{
    \includegraphics[width=200pt]{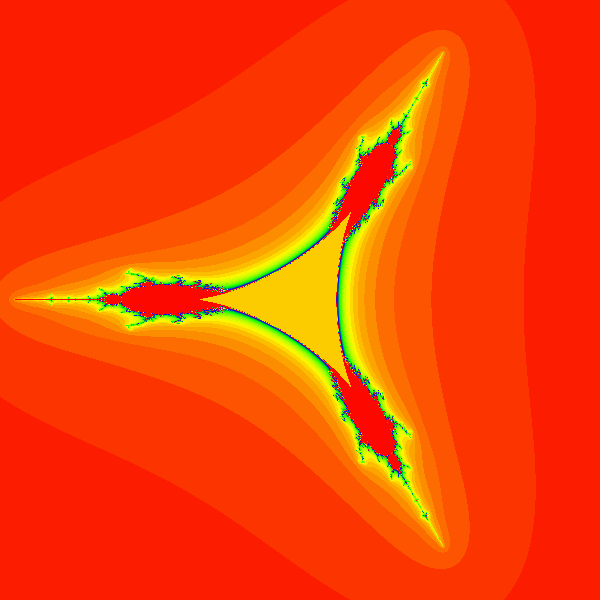}}
    \hspace{0.1in}
    \subfigure[\scriptsize{A swapping region.}  ]{
    \includegraphics[width=200pt]{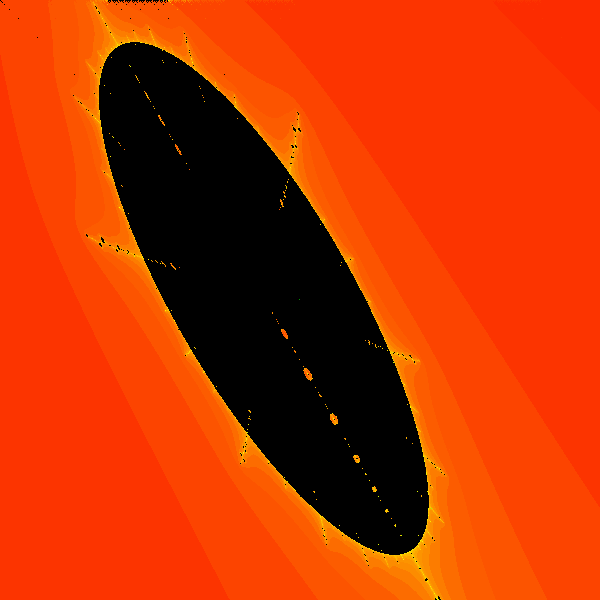}}
    \hspace{0.1in}
    \subfigure[\scriptsize{A copy of the Tricorn.}  ]{
    \includegraphics[width=200pt]{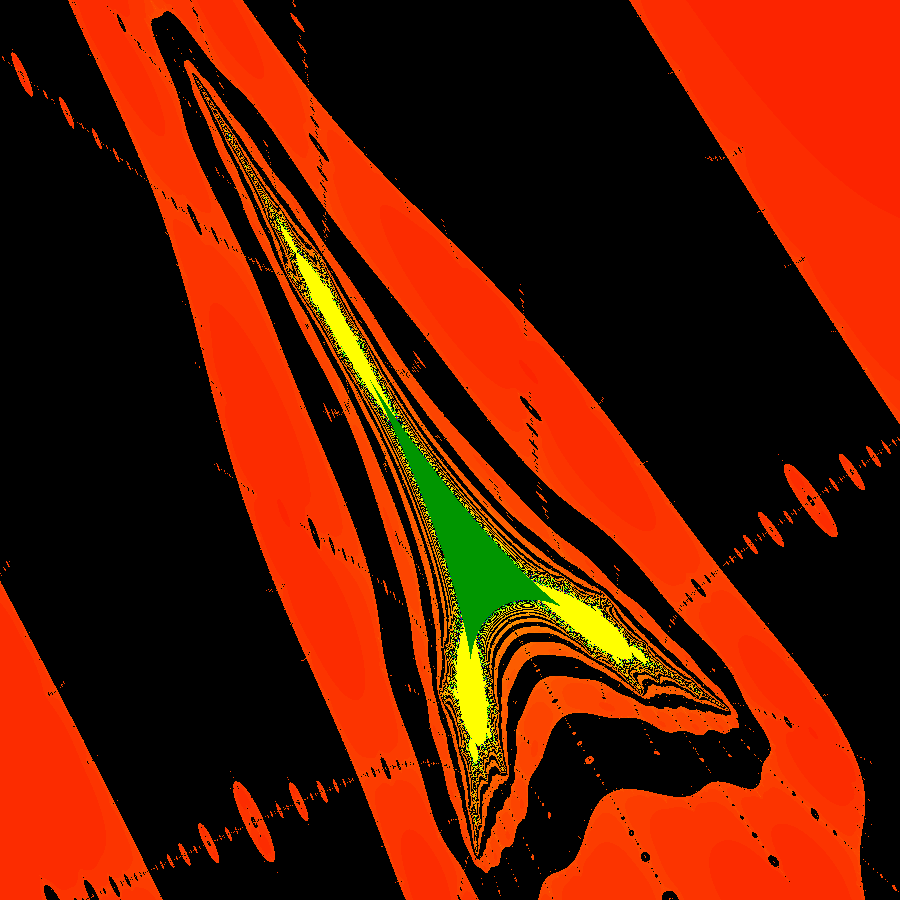}}
	\hspace{0.1in}
    \subfigure[\scriptsize{A copy of the Mandelbrot set.}  ]{
    \includegraphics[width=200pt]{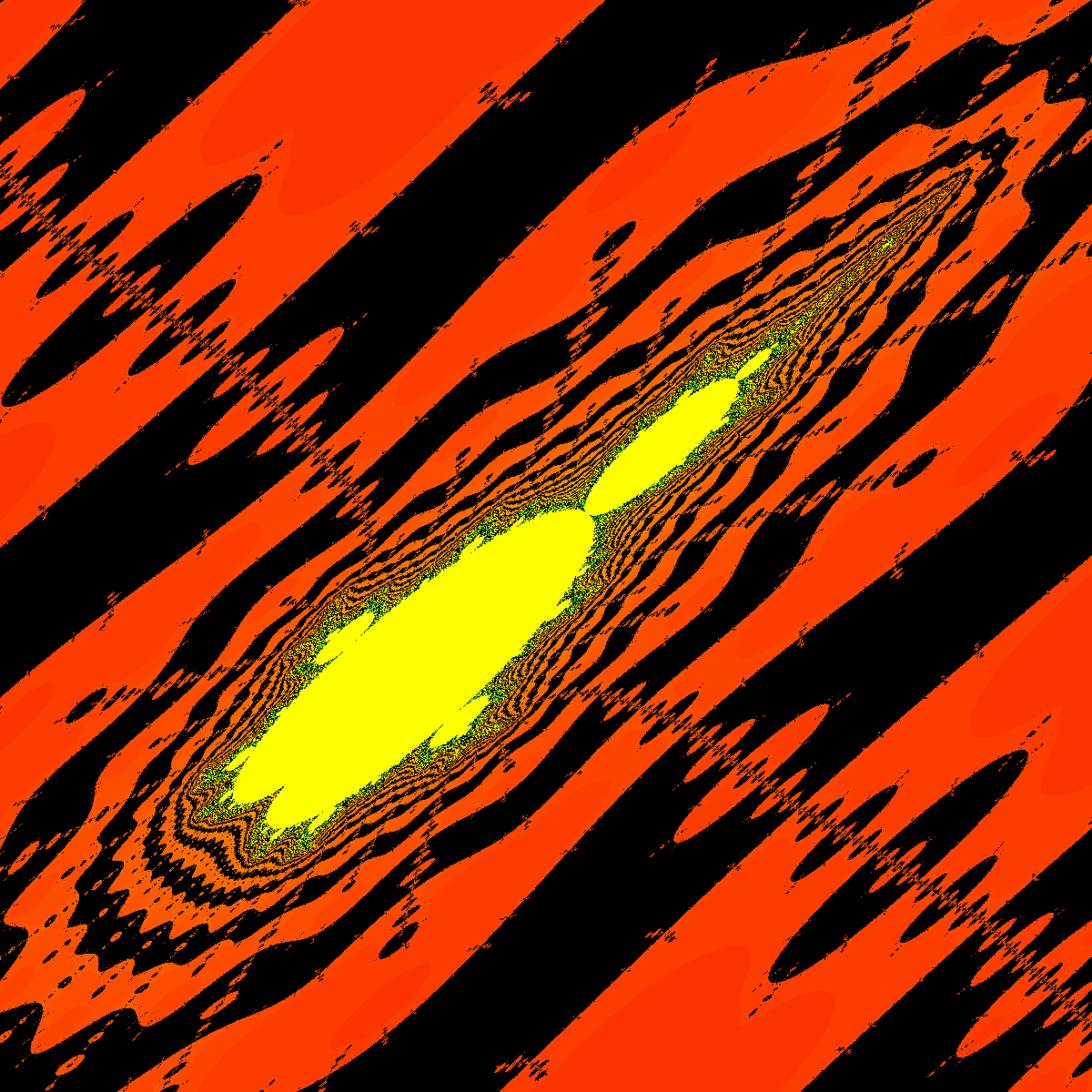}}
   
    \caption{\small{Figure (a) shows the Tricorn, the parameter plane of the antipolynomial $p_c(z)=\overline{z}^2+c$. Figure (b) shows the biggest  swapping region located in the tip of one of the three horns of the parameter space of the Blaschke family for $a\in$ $\left(-3.39603, -3.05761\right)\times\left(5.45471, 5.79312\right)$. It corresponds to the parameters bounded by the big black component. Figure (c) shows a zoom of (b) for which a copy of the Tricorn can be observed ($a\in(-3.22295, -3.22249)\times(5.58172, 5.58218)$). Figure (d) shows a copy of the Mandelbrot set inside another swapping component for $a\in(2.080306,2.080311)\times(1.9339165,1.9339215)$. In Figures (b), (c) and (d) red points correspond to parameters for which $\mathcal{O}(c_+)\rightarrow\infty$ whereas black points correspond to parameters for which $\mathcal{O}(c_+)\rightarrow 0$.  Green points correspond to bitransitive parameters (see Figure \ref{bitransdisjointswap} (left)), whereas yellow points correspond to disjoint parameters (see Figure \ref{bitransdisjointswap} (right)). The red and black annuli surrounding the copies of the Tricorn and the Mandelbrot set correspond to parameters for which $c_+$ enters and exists the unit disk a finite number of times before being captured by infinity or zero.}} 
    \label{swapregion}
\end{figure}

The goal of this section is to describe the dynamics for parameters which belong to swapping regions, i.e.\  parameters for which the critical orbit of $c_+$, $\mathcal{O}(c_+)$, enters the unit disk at least once (see Figure \ref{swapregion} (b)). Exploratory work shows that  small copies of the Tricorn, the connectedness locus of the antipolynomials $p_c(z)=\overline{z}^2+c$ (see \cite{Cr}, \cite{NaSh} and Figure \ref{swapregion} (a)), as well as small copies of the Mandelbrot set seem to appear embedded inside these regions (see Figures \ref{swapregion} (c) and (d)). They appear as the accumulation set of parameters for which $\mathcal{O}(c_+)$ enters and exits the unit disk more and more times. In the limit we may have parameters having attracting cycles which enter and exit the unit disk (see Figure \ref{swapregion} (c) and (d)). In this situation, we build, in Theorem~B, a polynomial like of degree either $2$ or $4$. Milnor \cite{Mi3}  worked with cubic polynomials with real parameters and described a similar structure than the one that we explain in Lemma \ref{periodehyp}. Using antipolynomial-like mappings (see Section \ref{santipol}) he related the dynamics of the bitransitive parameters with the ones of the Tricorn. We use his ideas to prove that when the polynomial-like map build in Theorem~B is of degree $4$, it is hybrid equivalent to a polynomial of the form $p^2_c(z)=\left(\overline{\overline{z}^2+c}\right)^2+c=\left(z^2+\overline{c}\right)^2+c$.

\clearpage

The following lemma tells us that swapping regions are disjoint from tongues (see Defintion \ref{deftongue}).

\begin{lemma}\label{notongueswapping}
A parameter $a$ with $|a|>2$ such that $B_a$ has an attracting or parabolic cycle in $\cercle$ cannot be swapping.
\end{lemma}
\proof
Assume that $B_a$ has an attracting cycle in $\cercle$. Let $\mathcal{A}$ be the maximal domain of the K{\oe}nigs linearization of the cycle (see \cite{Mi1}, Theorem 8.2 and Lemma 8.5). By symmetry, $c_{\pm}\in \partial\mathcal{A}$. By injectivity of the linearizer, $\mathcal{A}\cap\dis$ is mapped into $\mathcal{A}\cap\dis$ under $B_a$. Therefore, $\mathcal{O}(c_-)$ cannot exit the unit disk and the parameter is not swapping. The parabolic case is derived similarly taking $\mathcal{P}$ to be the maximal petal having the critical points on its boundary (see \cite{Mi1}, Theorems 10.9 and 10.15).
\endproof

We are interested in the hyperbolic regions contained in the swapping areas. We shall study them using the theories of polynomial and antipolynomial-like mappings. We focus in the bitransitive parameters which, for $|a|>2$, are necessarily inside swapping regions (see Theorem \ref{clashyp} (b)).

We recall some notation from Section \ref{introblas}. For $|a|>2$ the unit circle has two preimages different from itself and not intersecting $\cercle$, say $\gamma_i\subset\dis$ and $\gamma_e\subset \com\setminus\overline{\dis}$ (see Figure~\ref{esquemapunts}~(d)). The map $B_a$ sends $\gamma_i$ and $\gamma_e$ bijectively to $\cercle$. Moreover, let $\Omega_i$ be the region bounded by $\gamma_i$ and contained in $\dis$ and let $\Omega_e$ be the region bounded by $\gamma_e$ and contained in $\com\setminus\overline{\dis}$. The maps $B_a|_{\Omega_i}:\Omega_i\rightarrow \widehat{\com}\setminus\overline{\dis}$ and  $B_a|_{\Omega_e}:\Omega_e\rightarrow \dis$ are conformal (there is only one preimage of $z=\infty$ in $\dis$ and one preimage of $z=0$ in $\com\setminus\overline{\dis}$), so $c_{\pm}\notin \Omega_{i}\cup \Omega_e$.

We now prove a lemma which deals with the period of attracting and parabolic cycles for parameters inside swapping regions.

\begin{lemma}\label{periodehyp}
Let $a$, $|a|>2$, be a parameter inside a swapping region. If $B_a$ has an attracting or parabolic cycle, then its period is at least 3. Moreover, if $a$ is bitransitive, the period is even.
\end{lemma}

\proof
First of all notice that, from Lemma \ref{notongueswapping} and invariance of $\cercle$, no component of the basin of attraction of the cycle can intersect neither $\gamma_i$ nor $\gamma_e$. A parabolic or attracting cycle needs to have a critical point in its immediate basin of attraction. The component in which the critical point lies is contained neither in $\Omega_i$ nor in $\Omega_e$. Moreover, since the periodic cycle  needs to enter and exit the unit disk, the immediate basin of attraction of the cycle has a component in $\Omega_i$ and another one in $\Omega_e$. Then, the immediate basin of attraction has at least three different components and, hence, the cycle has at least period three.

Now assume that $a$ is bitransitive. Suppose without lose of generality that the component which contains $c_+$ is mapped under $k>0$ iterates to the component which contains $c_-$. Because of  symmetry, the first return map from the component of $c_-$ to component of $c_+$ also takes $k$ iterates. Hence, the period of the attracting cycle is $2k$.
\endproof

We proceed now to prove Theorem B.

\begin{proof}[Poof of Theorem B]

First of all notice that, due to Lemma \ref{notongueswapping}, $A^*(<z_0>)$  neither intersects $\gamma_e$ nor $\gamma_i$. Since the cycle enters and exits the unit disk, $A^*(<z_0>)$ has at least one connected component entirely contained in $\Omega_e$. Let $A^*(z_0)$ be the connected component of  $A^*(<z_0>)$ containing $c_+$. Let $n_0\in \nat$ be minimal such that $B_{a_0}^{n_0}(z_0)=z_{n_0}\in\Omega_e$. Let $\mathcal{S}_0$  be the connected component of $B_{a_0}^{-n_0}(\Omega_e)$, containing $c_+$ (and hence $A^*(z_0)$). The set $\mathcal{S}_0$ is simply connected by Corollary \ref{connpreim} since $\Omega_e$ is simply connected and $B_{a_0}^{n_0}|_{\mathcal{S}_0}$ has a unique critical point. Recursively define $\mathcal{S}_n$ to be the connected component of $B_{a_0}^{-1}(\mathcal{S}_{n-1})$ containing the point $z_{-n}$ of the cycle (recall that the subindexes of the cycle are taken in $\mathbb{Z}/p\mathbb{Z}$). Again by Corollary \ref{connpreim}, the components $\mathcal{S}_n$ are simply connected for all $n>0$. Let $p_0\in \nat$ be the minimal  such that $c_+\in\mathcal{S}_{p_0}$. Since $\mathcal{S}_{p_0-n_0}\subset\Omega_e$, we have that $\mathcal{S}_{p_0}\subset\mathcal{S}_0$. Notice that $p_0$ is a divisor of $p$.

 The map $B_{a_0}|_{\mathcal{S}_n}:\mathcal{S}_n\rightarrow\mathcal{S}_{n-1}$ is conformal if $\mathcal{S}_n$ contains no critical point and 2-to-1 if it contains $c_+$ or $c_-$ (it cannot contain both critical points at the same time since $\mathcal{S}_n\cap\cercle=\emptyset$). Hence, the triple $(B_{a_0}^{p_0};\mathcal{S}_{p_0},\mathcal{S}_0)$ is a polynomial-like map of degree $4$ or $2$ depending on whether there is some $\mathcal{S}_{q_0}$ containing $c_-$ or not. As in  Lemma \ref{periodehyp}, if such $q_0$ exists, $p_0=2q_0$. Notice that, if the parameter is bitransitive, this $q_0$ exists and, therefore, the degree is $4$. Since the condition $c_+(a)\in\mathcal{S}_{p_0}(a)\subset\mathcal{S}_0(a)$ is open, the polynomial-like map can be defined for an open set of parameters $W$ around $a_0$. From now on we consider $a\in W$.

\begin{figure}[hbt!]
    \centering

\def\svgwidth{250pt}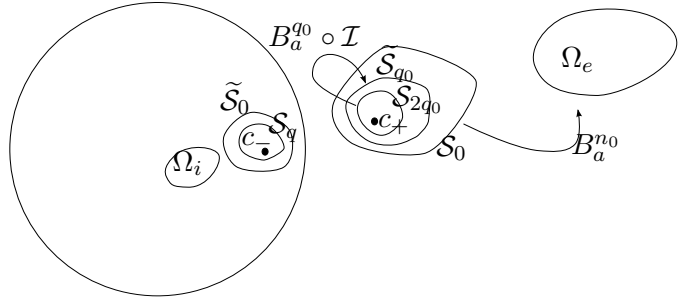  
    \caption{\small{ Sketch of the situation described in Theorem B for the degree 4 case.}}
    \label{antipol2}
\end{figure}

 We now use antipolynomial-like mappings to see that, in the case of degree $4$ polynom\-ial-like mapping, the degree $4$ polynomial to which $(B_{a}^{2q_0};\mathcal{S}_{2q_0},\mathcal{S}_0)$ is hybrid equivalent can be taken of the form $p_c^2(z)=\left(z^2+\overline{c}\right)^2+c$. See Section \ref{santipol} for an introduction to antipolynomial-like mappings. We proceed to construct an antypolinomial-like map $(f;\mathcal{S}_{2q_0}, \mathcal{I}(\mathcal{S}_{q_0}))$ of degree 2, where $\mathcal{I}(z)=1/\overline{z}$ denotes the reflection with respect to $\cercle$. This antipolynomial-like map is hybrid equivalent to an  antipolynomial of the form  $\overline{z}^2+c$. The result then follows if $f(\mathcal{I}(\mathcal{S}_{q_0}))= \mathcal{S}_0$ and $(f^2;\mathcal{S}_{2q_0}, \mathcal{S}_0)=(B_a^{2q_0};\mathcal{S}_{2q_0}, \mathcal{S}_0)$.

Define $\widetilde{\mathcal{S}}_q=\mathcal{I}(\mathcal{S}_q)$, where $q\in\nat$. It is easy to see that $\overline{\mathcal{S}_{q_0}}\subset\widetilde{\mathcal{S}}_0$. Indeed, taking $n_0$ as in the definition of $\mathcal{S}_0$,  by symmetry, $B_a^{n_0}(\widetilde{\mathcal{S}}_0)=\Omega_i$. Since $B_a^{n_0}(\overline{\mathcal{S}_{q_0}})$ is contained in $\Omega_i$, we conclude that $\overline{\mathcal{S}}_{q_0}\subset \widetilde{\mathcal{S}}_0$ (see Figure \ref{antipol2}). From $\overline{\mathcal{S}_{q_0}}\subset\widetilde{\mathcal{S}}_0$ we can deduce that $\overline{\mathcal{S}_{2q_0}}\subset\widetilde{A}_{q_0}$. Finally, take $f= \mathcal{I}\circ B_a^{q_0}$. Since $B_a=\mathcal{I}\circ B_a\circ \mathcal{I}$ we have that $f^2=B_a^{2q_0}$. Then, the antipolynomial-like map $(\mathcal{I}\circ B_a^{q_0};\mathcal{S}_{2q_0}, \widetilde{\mathcal{S}}_{q_0}=\mathcal{I}(\mathcal{S}_{q_0}))$ satisfies the desired conditions.

\end{proof}

\begin{figure}[hbt!]
    \centering

    \subfigure{
    \includegraphics[width=190pt]{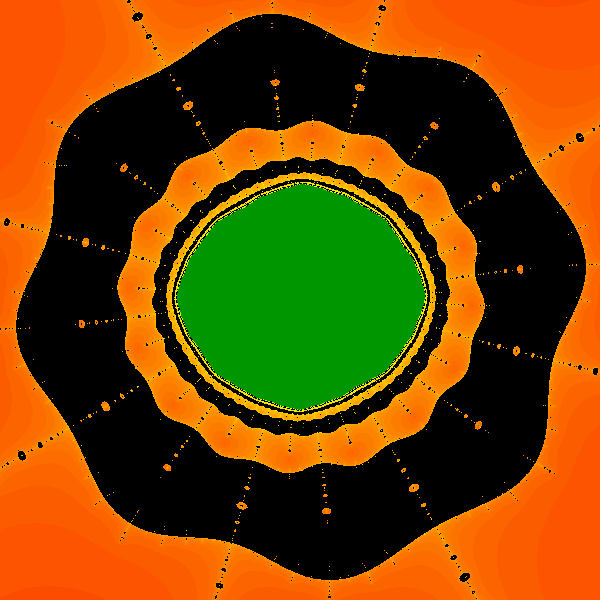}}
     \hspace{0.1in}
     \subfigure{
 	\includegraphics[width=190pt]{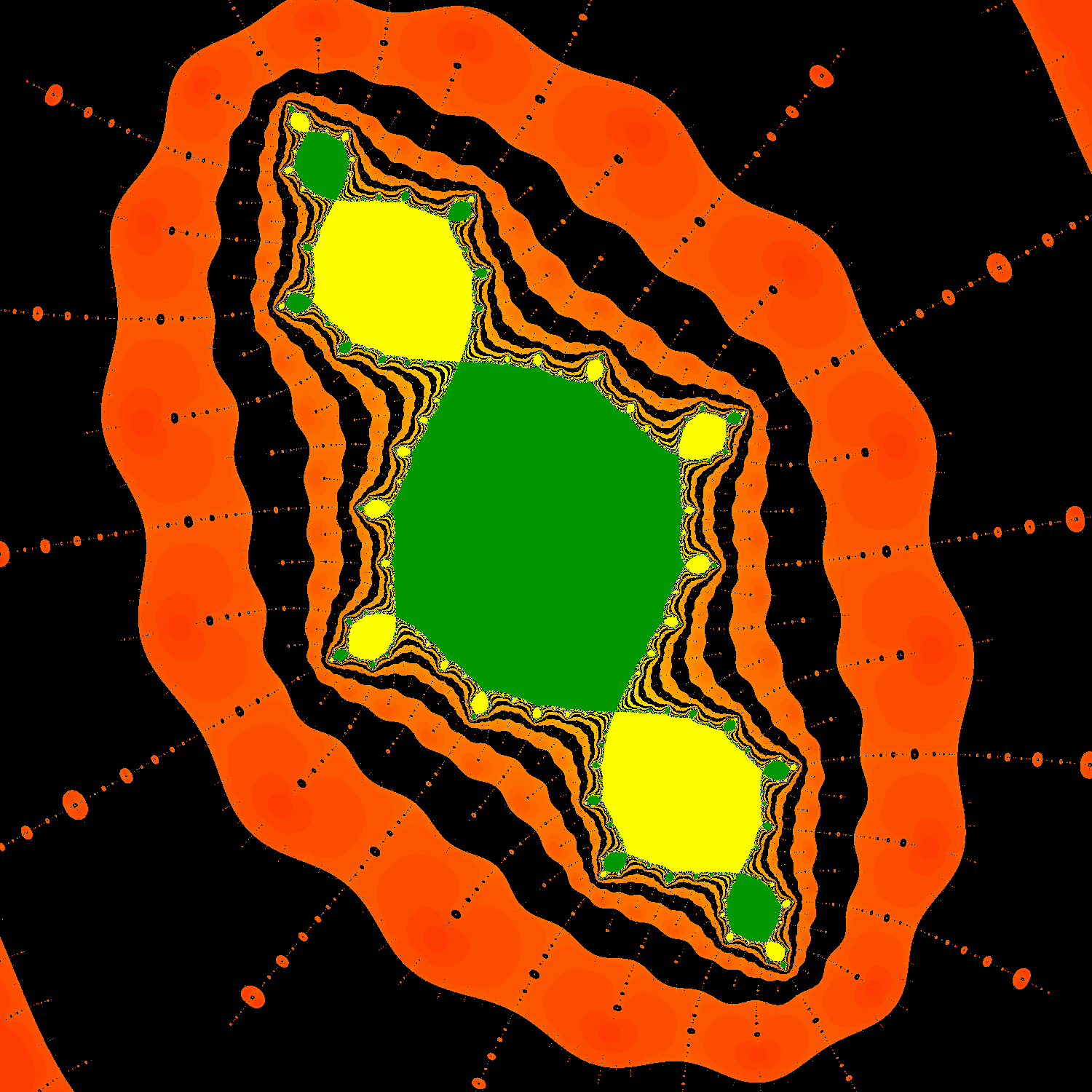}}
   
    \caption{\small{The left figure shows a connected component of the bitansitive cycle of $B_{a_1}$, where $a_1=-3.22271+5.58189i$. The right figure shows the  dynamical plane of $B_{a_2}$, where $a_2=-3.22278+5.58202i$ is a disjoint swapping parameter. Colours work as in Figure \ref{1<a<2}. Notice that, surrounding the basin of attraction, appear some black and red annuli. These annuli correspond to points which enter and exit $\dis$ a finite number of times and then are attracted to zero or infinity.}} 
    \label{bitransdisjointswap}
\end{figure}

Theorem B tells us that all bitransitive parameters contained in swapping regions can be related to the dynamics of $p_c^2(z)$, where $p_c(z)=\overline{z}^2+c$, since the polynomial-like map has degree $4$. However, notice that if an antipolynomial $p_c(z)$ has an attracting cycle of even period $2q$, then $p_c^2(z)$ has two disjoint attracting cycles of period $q$. Therefore,  there are also parameters  which appear as disjoint parameters in the copies of the Tricorn in the parameter plane of the Blaschke family (see Figure \ref{swapregion} (c) and Figure \ref{bitransdisjointswap} (right)). These disjoint parameters also lead to degree $4$ polynomial-like maps. Finally, the polynomial-like maps of degree $2$ obtained from the other disjoint parameters are hybrid equivalent to quadratic polynomials $z^2+c$. These parameters correspond to the ones inside the small copies of the Mandelbrot set observed by means of numerical computations (see Figure \ref{swapregion} (d)).

\subsection{Parametrizing hyperbolic components: Proof of Theorem C}\label{parametrizehyp}

The aim of this section is to study the multiplier map of the bitransitive and disjoint hyperbolic components of the Blaschke family $B_a$  for parameters $a$ such that $|a|>2$. Recall that a hyperbolic component is a connected component of the set of parameters for which $B_a$ is hyperbolic. The section is structured as follows: first we prove a
proposition that is useful later on, then we see that the multiplier map is a homeomorphism between any disjoint hyperbolic component and the unit disk proving Theorem C and finally we study the bitransitive case.

The following proposition tells us that, given $B_a$ with $|a|>2$, the boundaries of the connected components of the basin of attraction of every attracting cycle contained in $\com^*\setminus\cercle$ are Jordan curves. The result is a direct consequence of the relation of the family $B_a$ with polynomials which has been described in Proposition \ref{surgerymilnor} and Theorem B, respectively.

\begin{propo}\label{jordancurveattr}
Assume that $B_{a}$  has an attracting cycle $<z_0>$ which is contained in $\com^*\setminus\cercle$. Then, the boundaries of the connected components of the basin of attraction $A(<z_0>)$ are Jordan curves.
\end{propo}
\proof

It follows from the hypothesis of the proposition that $|a|>2$ since for $1<|a|\leq 2$ any attracting cycle other than $z=0$ or $z=\infty$ is contained in $\cercle$ and for $|a|<1$ there are no attracting cycles in $\com^*$.  It follows from Proposition~\ref{surgerymilnor} and Theorem~B  that the closure of every connected component of $A^*(<z_0>)$ is homeomorphic to the closure of a connected component of a bounded attracting cycle of a polynomial. Since the boundary of every bounded Fatou component of a polynomial other than a Siegel disk is a Jordan curve (see \cite{RoYi}), the boundary of every connected component of $A^*(<z_0>)$ is also a Jordan curve. Finally, since all critical points are contained in the immediate basins of attraction of attracting cycles, the closure of every connected component $U$ of $A(<z_0>)\setminus A^*(<z_0>)$ is mapped homeomorphically to the closure of a connected component of $A^*(<z_0>)$ and, therefore, $\partial U$ is a Jordan curve too.

\endproof

We now prove Theorem~C. The main idea of the proof is to build a local inverse of the multiplier map $\Lambda$ around every $\Lambda(a_0)\in\dis$. It is done performing a cut and paste surgery to  change the multiplier of the attracting cycle using a degree 2 Blaschke product with an attracting cycle of the desired multiplier as a model (see \cite{BF}, Chapter 4.2).

\proof[Proof of Theorem C]

Consider the family of degree 2 Blaschke products 

$$b_{\lambda}(z)=z\frac{z+\lambda}{1+\overline{\lambda}z},$$

\noindent where $\lambda\in\dis$. They have $0$ and $\infty$ as attracting fixed points of multipliers $\lambda$ and $\overline{\lambda}$, respectively. The only other fixed point $\frac{1-\lambda}{1-\overline{\lambda}}\in\cercle$ is repelling. The multiplier $\lambda$ and the repelling fixed point in $\cercle$ determine univalently the map $b_{\lambda}$ since any holomorphic self-map of degree 2 of $\dis$ has the previous form. Its Julia set satisfies $\mathcal{J}(b_{\lambda})=\cercle$. Furthermore, if $|\lambda|<r<1$, $D_{\lambda}=b_{\lambda}^{-1}(\dis_r)$ is a simply connected open set which compactly contains the disk of radius $r$, $\dis_r=\{z, |z|<r\}$, whereas $b_{\lambda}(\dis_r)$ is compactly contained in $\dis_r$.

Let $a_0\in U$ and let $\lambda_0$ be the multiplier of the attracting cycle $<z_0>$ of period $p$ of $B_{a_0}$ such that $c_+\in A^*(z_0)$. Since there is no other critical point in $A^*(<z_0>)\setminus A^*(z_0)$ and $\partial A^*(z_0)$ is a Jordan curve (Proposition \ref{jordancurveattr}), the map $B_{a_0}^p: \overline{A^*(z_0)}\rightarrow \overline{A^*(z_0)}$ has degree 2 and a unique fixed point $z'_0$ in $\partial A^*(z_0)$. Let $\mathcal{R}:A^*(z_0)\rightarrow \dis$ be the Riemann map sending $z_0$ to $0$ and $z'_0$ to $\frac{1-\lambda_0}{1-\overline{\lambda_0}}$. The map  $\mathcal{R}\circ B_{a_0}^p \circ \mathcal{R}^{-1}$ is, by construction, the restriction to $\dis$ of the Blaschke product $b_{\lambda_0}(z)$. Fix $r'$ and $r$ so that $0\leq |\lambda_0|< r'<r<1$. We proceed now to perform a surgery to the product $b_{\lambda_0}$ which  changes the multiplier of the attracting fixed point $0$ to $\lambda$ for any $|\lambda|<r'$.

Let $D_{\lambda_0}=b_{\lambda_0}^{-1}(\dis_r)$ and let $A_{\lambda_0}$ denote the annulus $D_{\lambda_0}\setminus \dis_r$. Define $g_{\lambda}:\dis\rightarrow\dis$ as

$$g_{\lambda}=\left\{\begin{array}{lcl}
b_{\lambda_0} &  \mbox{ on } & \dis\setminus D_{\lambda_0}\\
b_{\lambda} &  \mbox{ on } & \overline{\dis}_r\\
h_{\lambda} &  \mbox{ on } &A_{\lambda_0},
\end{array}\right. $$

\noindent where $h_{\lambda}$ is chosen to be a quasiconformal map which interpolates $b_{\lambda}$ and $b_{\lambda_0}$ depending continuously on $\lambda$. Such a interpolating map can be taken since the boundary maps $g_{\lambda}|_{\partial A_{\lambda_0}}$ are degree 2 analytic maps on analytic curves. The inner boundary map  depends continuously on $\lambda$ whereas the outer map is independent of it. Therefore, the map $h_{\lambda}:A_{\lambda_0}\rightarrow A_{\lambda}$, where $A_{\lambda}$ denotes the annulus  $\dis_r\setminus b_{\lambda}(\dis_r)$, can be chosen to be a quasiconformal covering map of degree 2 which depends continuously on $\lambda$ (see Exercise 2.3.3 of \cite{BF}). We define recursively a $g_{\lambda}$-invariant almost complex structure $\tilde{\sigma}_{\lambda}$ as

$$\widetilde{\sigma}_{\lambda}=\left\{\begin{array}{lcl}
\sigma_0 &  \mbox{ on } & \dis_r\\
h_{\lambda}^*\sigma_0 &  \mbox{ on } & A_{\lambda_0}\\
(b_{\lambda_0}^n)^*\widetilde{\sigma}_{\lambda} &  \mbox{ on } &b_{\lambda_0}^{-n}(A_{\lambda_0}),
\end{array}\right.$$

\noindent where $\sigma_0$ denotes the standard complex structure. Notice that, since any $z\in\dis$ can go at most once trough $A_{\lambda_0}$, $\widetilde{\sigma}_{\lambda}$ has bounded dilatation. Indeed, $||h_{\lambda}^*\sigma_0||_{\infty}:=k(\lambda)<1$ since $h_{\lambda}$ is quasiconformal and the pull backs $(b_{\lambda_0}^n)^*$ do not increase the dilatation. Moreover, $\lambda\rightarrow\widetilde{\sigma}_{\lambda}(z)$ varies continuously with $\lambda$ for all $z\in \dis$. Notice also that the almost complex structures have dilatation uniformly bounded by $k:=max_{|\lambda|\leq r'}k(\lambda)<1$ for all $\lambda\in \dis_{r'}$.

Once we have performed the multiplier surgery in the degree 2 Blaschke model, we glue it in $B_{a_0}$. This is done preserving the symmetry of the family, i.e.\ the new map is preserved under pre and post composition by $\mathcal{I}(z)=1/\overline{z}$. Define the model map $F_{\lambda}$ as

$$F_{\lambda}=\left\{\begin{array}{lcl}
(B_{a_0}^{-1})^{(p-1)}\circ\mathcal{R}^{-1}\circ g_{\lambda}\circ \mathcal{R} &  \mbox{ on } & A^*(z_0)\\
\mathcal{I}\circ F_{\lambda} \circ \mathcal{I} &  \mbox{ on } & \mathcal{I}(A^*(z_0))\\
B_{\lambda_0} &  \mbox{ elsewhere,} &
\end{array}\right. $$

\noindent where $(B_{a_0}^{-1})^{(p-1)}$  denotes $\displaystyle B_{a_0}^{-1}\circ \mathop{\cdots}^{p-1}\circ B_{a_0}^{-1}$. It is well defined since $B_{a_0}:A^*(z_i)\rightarrow A^*(z_{i+1})$ is conformal for every $i\neq0$. The map $F_{\lambda}$ depends continuously on $\lambda$, is symmetric with respect to $\cercle$ and holomorphic everywhere except in $\mathcal{R}^{-1}(A_{\lambda_0})\cup \mathcal{I}(\mathcal{R}^{-1}(A_{\lambda_0}))$. Notice also that the periodic cycle $<z_0>$ of $F_{\lambda}$ has multiplier $\lambda$ and that $F_{\lambda}$ has a unique critical point in $\com\setminus\dis$ which depends continuously on $\lambda$. We define recursively an $F_{\lambda}$-invariant almost complex structure $\sigma_{\lambda}$ as

$$\sigma_{\lambda}=\left\{\begin{array}{lcl}
\mathcal{R}^*\widetilde{\sigma}_{\lambda} &  \mbox{ on } & A^*(z_0)\\
\mathcal{I}^*\sigma_{\lambda} &  \mbox{ on } & \mathcal{I}(A^*(z_0))\\
(B_{a_0}^n)^*\sigma_{\lambda} &  \mbox{ on } & B_{a_0}^{-n}(A^*(z_0))\setminus A^*(z_0)\\
(B_{a_0}^n)^*\sigma_{\lambda} &  \mbox{ on } & B_{a_0}^{-n}(\mathcal{I}(A^*(z_0)))\setminus \mathcal{I}(A^*(z_0))\\
\sigma_0 &  \mbox{ elsewhere.} &
\end{array}\right. $$

Notice that we are pulling back the almost complex structure $\sigma_{\lambda}$ by the antiholomorphic map $\mathcal{I}(z)$ (see Section 1.2.1 of \cite{BF} for an introduction to pull backs under orientation reversing maps) and that, since we only pull back under holomorphic and antiholomorphic maps, $||\widetilde{\sigma}_{\lambda}||_{\infty}=||\sigma_{\lambda}||_{\infty}$. By Theorem \ref{conjpe}, $B_{a_0}|_{\cercle}$ is conjugate to the doubling map and, therefore, has a unique fixed point $x_0\in\cercle$. Let $\phi_{\lambda}:\hat{\com}\rightarrow\hat{\com}$ be the integrating map of $\sigma_{\lambda}$ obtained from the Measurable Riemann Mapping Theorem (see \cite{Ah}) and normalized so that it fixes $0$, $x_0$ and $\infty$. Since $\sigma_{\lambda}(z)$ depends continuously on $\lambda$ for all $z\in\hat{\com}$ with dilatation which is uniformly bounded away from 1, the map $\phi_{\lambda}(z)$ depends continuously on $\lambda$ for each $z\in\hat{\com}$.  It follows from the uniqueness of the integrating map and the symmetry of $\sigma_{\lambda}$ with respect to $\cercle$ (i.e.\  $\sigma_{\lambda}=\mathcal{I}^*\sigma_{\lambda}$) that $\phi_{\lambda}$ is also symmetric with respect to $\cercle$. Therefore,  $\tilde{B}_{\lambda}=\phi_{\lambda} \circ F_{\lambda}\circ \phi_{\lambda}^{-1}$ is a degree 4 holomorphic map of $\hat{\com}$ symmetric with respect to $\cercle$ which has $z=0$ and $z=\infty$ as superattracting cycles of local degree $3$. Therefore, $\tilde{B}_{\lambda}$ is a Blaschke product of the form $B_{\tilde{a}(\lambda), \tilde{t}(\lambda)}$  (\ref{blasformula2}). Since $F_{\lambda}$ has a unique critical point in $\com\setminus\dis$ which depends continuously on $\lambda$ and $\phi_{\lambda}(z)$ depends continuously on $\lambda$ for each $z\in\hat{\com}$, $B_{\tilde{a}, \tilde{t}}$ has a unique critical point $c_+(\lambda)\in \com\setminus\dis$ which depends continuously on $\lambda$.  Therefore, since $B_{\tilde{a}(\lambda), \tilde{t}(\lambda)}$ fixes $x_0\in\cercle$, we have by Lemma \ref{continuitya} that $\tilde{a}(\lambda)$ and $\tilde{t}(\lambda)$ depend continuously on $\lambda$. Finally, by lemma \ref{conjblas}, $B_{\tilde{a}(\lambda), \tilde{t}(\lambda)}$ is conjugate to a Blaschke product $B_{a(\lambda)}$ (\ref{blasformula}), where $a(\lambda)=\tilde{a}(\lambda)e^{\frac{2\pi i \tilde{t}({\lambda})}{3}}$ depends continuously on $\lambda$.

To finish the proof we check that  $a(\lambda_0)=a_0$ and, therefore, every $a(\lambda)$ belongs to the same hyperbolic component $U$ as $a_0$. We have not justified that the  quasiconformal interpolating map $h_{\lambda_0}$ equals $b_{\lambda_0}$ and, hence, $B_{a_0}$  and $B_{\tilde{a}(\lambda_0),\tilde{t}_0}$ might be distinct. However, the integrating map $\phi_{\lambda_0}$ is a conformal conjugacy between them in $\widehat{\com}\setminus\overline{A(<z_0>)\cup A(< \mathcal{I}(z_0)>)}$ and is a quasiconformal conjugacy in a neighbourhood of their Julia sets. Define $\widetilde{\phi}_{\lambda_0}$ to be the conformal map from $A(<z_0>)\cup A(<\mathcal{I}(z_0)>)$ to $A(<\phi_{\lambda_0}(z_0)>)\cup A(<\mathcal{I}(\phi_{\lambda_0}(z_0))>)$ such that, restricted to every connected component,  coincides with the Riemann map normalized so that the attracting cycles $<z_0>$,  $<\mathcal{I}(z_0)>$, $<z'_0>$ and $<\mathcal{I}(z'_0)>$ are mapped to $< \phi_{\lambda_0}(z_0)>$, $<\mathcal{I}( \phi_{\lambda_0}(z_0))>$, $< \phi_{\lambda_0}(z'_0)>$ and $<\mathcal{I}( \phi_{\lambda_0}(z'_0))>$ and their preimages are in correspondence. Since $B^p_{a_0}$  is conjugate to $b_{\lambda_0}$ (resp.\ $b_{\overline{\lambda_0}}$) in $A^*(z_0)$ (resp.\ $A^*(\mathcal{I}(z_0))$)  and so is $B^p_{\tilde{a}(\lambda_0),\tilde{t}_0}$ in $A^*(\phi_{\lambda_0}(z_0))$ (resp.\ $A^*(\mathcal{I}(\phi_{\lambda_0}(z_0)))$), the conformal map $\widetilde{\phi}_{\lambda_0}$ is a conjugacy.  Moreover, it extends to the boundary of every connected component of the basins of attraction since they are Jordan domains by Proposition \ref{jordancurveattr}. Given that $\phi_{\lambda_0}$ and $\widetilde{\phi}_{\lambda_0}$ conjugate $B_{a_0}$ and $B_{\tilde{a}(\lambda_0),\tilde{t}_0}$ in $\partial A(<z_0>)\cup \partial A(<\mathcal{I}(z_0)>)\subset \mathcal{J}(B_{a_0})$ they coincide since they map periodic points to periodic points. Consequently, the map $\varphi_{\lambda_0}$ defined as $\phi_{\lambda_0}$ in $\hat{\com}\setminus\left( A(<z_0>)\cup  A(<\mathcal{I}(z_0)>)\right)$ and $\widetilde{\phi}_{\lambda_0}$ in $ A(<z_0>)\cup A(<\mathcal{I}(z_0)>)$ is a global conjugacy. Moreover, since $\varphi_{\lambda_0}$ is quasiconformal in $\widehat{\com}\setminus\left(\partial A(<z_0>)\cup\partial A(<\mathcal{I}(z_0)>)\right)$, coincides with $\phi_{\lambda_0}$ in  $\mathcal{J}(B_{\lambda_0})$ and $\phi_{\lambda_0}$ is quasiconformal in a neighbourhood of $\mathcal{J}(B_{\lambda_0})$,  $\varphi_{\lambda_0}$ is quasiconformal by Rickman's Lemma (cf.\ \cite{DH1}, \cite{BF}).  Since $\varphi_{\lambda_0}$ is conformal a.e.\ in $\widehat{\com}$, it is 1-quasiconformal and therefore conformal map of $\widehat{\com}$ by Weyl's Lemma. Since $\varphi_{\lambda_0}$ fixes $0$ and $\infty$, leaves $\cercle$ invariant and fixes $x_0\in\cercle$, we conclude that $\varphi_{\lambda_0}$ is the identity and $B_{\tilde{a}(\lambda_0),\tilde{t}_0}=B_{a_0}$.

For every $a_0\in U$ we have constructed a continuous local inverse to the multiplier map $\Lambda:U\rightarrow\dis$. Therefore, $\Lambda$ is a homeomorphism.
\endproof

We finish this section giving some ideas of what happens with bitransitive parameters (see Figure \ref{swapregion} (c)). It follows from Theorem~B that these parameters are strongly related to the quadratic antiholomorphic polynomials $p_c(z)=\overline{z}^2+c$. Indeed, the polynomial-like map constructed in Theorem~B is hybrid equivalent to a degree $4$ polynomial of the form $p_c^2(z)$ with a bitransitive attracting cycle. Therefore, the polynomial $p_c(z)$ also has an attracting cycle of odd period since, otherwise, the attracting $p_c^2(z)$ would have two disjoint attracting cycles.
Nakane and Schleicher $\cite{NaSh}$ studied the parameter plane of the antipolynomials $p_{c,d}=\overline{z}^d+c$ and, in particular, $p_{c,2}(z)=p_c(z)$.  If the period of the cycles of a hyperbolic component was even, they proved a result analogous to Theorem~C. They also showed that the multiplier map is not a good model for the odd period hyperbolic components. The reason why the multiplier map is not good for this case is the fact that the antiholomorphic multiplier $\frac{\partial}{\partial\overline{z}}f^k(z_0)$ of a cycle $<z_0>$ of odd period $k$ of an antiholomorphic map $f(z)$ is not a conformal invariant, only its absolute value is. They proved that the multiplier of the period $k$ cycle $<z_0>$ of the holomorphic map $f^2(z)$ equals the square of the absolute value of the previous antiholomorphic multiplier. Given a bitrinsitive hyperbolic component $U$ of $p_c^2(z)$, it also follows from their work that the set of parameters $c\in U$ for which the attracting cycle has multiplier $\lambda\in(0,1)$ is a Jordan curve and that $U$ contains a unique parameter $c_0$ for which the cycle is superattracting. We expect a similar result for bitransitive hyperbolic components of the Blaschke family $B_a$, but we only prove, for the sake of completeness, the following result.


\begin{propo}\label{bitransmult}
Let $<z_0>$ be a bitransitive cycle of a Blaschke product $B_a$ as in (\ref{blasformula}) with $|a|>2$. Then it has non-negative  real multiplier.
\end{propo}
\proof
By Lemma \ref{periodehyp}, the cycle $<z_0>$ has even period $2q$. Let $\mathcal{I}(z)=1/\overline{z}$. By symmetry, $\mathcal{I}(B^q_a(z_0))=z_0$. Therefore, $z_0$ is a fixed point of the antiholomorphic rational map $f=\mathcal{I}\circ B^q_a$. Moreover, $B_a^{2q}=f^2$. Therefore, the multiplier of the cycle is given by
$$\frac{\partial}{\partial z} B_a^{2q}(z_0)=\frac{\partial}{\partial\overline{z}}f(z_0)\cdot \overline{\frac{\partial}{\partial\overline{z}}f(z_0)}=\left|\frac{\partial}{\partial\overline{z}}f(z_0)\right|^2.$$

\endproof

\bibliography{tesi}
\bibliographystyle{amsalpha}

\end{document}